\documentclass{mcom-l}
\usepackage{caption}
\usepackage{subcaption}
\usepackage{enumitem}
\usepackage{amsmath}
\usepackage{amssymb}
\usepackage{bbm}
\usepackage{verbatim}
\usepackage{graphicx,epstopdf}
\usepackage{array}
\usepackage{tikz-cd}
\usetikzlibrary{arrows}
\tikzset{
   commutative diagrams/.cd,
   arrow style=tikz,
   diagrams={>=latex},
   every label/.append style = {font = \normalsize}}
\usepackage{xfrac}
\usepackage{faktor}
\newcolumntype{R}{>{$}r<{$}} %
\newcolumntype{V}[1]{>{[\;}*{#1}{R@{\;\;}}R<{\;]}} %
\usepackage{braket,amsfonts}
\usepackage{multirow}
\usepackage{booktabs}
\usepackage{graphbox}
\usepackage{overpic}
\usepackage{enumitem}
\usepackage{url}

\usepackage{pgf,tikz,pgfplots}
\pgfplotsset{compat=1.15}
\usetikzlibrary{shapes,arrows}
\usetikzlibrary{calc,patterns,angles,quotes}

\usepackage{algorithm}

\usepackage{algorithmic}

\usepackage{graphicx,epstopdf}

\usepackage{amsopn}

\usepackage{xspace}
\usepackage{bold-extra}
\usepackage[most]{tcolorbox}

\colorlet{texcscolor}{blue!50!black}
\colorlet{texemcolor}{red!70!black}
\colorlet{texpreamble}{red!70!black}
\colorlet{codebackground}{black!25!white!25}

\tcbset{%
  colframe=black!75!white!75,
  coltitle=white,
  colback=codebackground, 
  colbacklower=white, 
  fonttitle=\bfseries,
  arc=0pt,outer arc=0pt,
  top=1pt,bottom=1pt,left=1mm,right=1mm,middle=1mm,boxsep=1mm,
  leftrule=0.3mm,rightrule=0.3mm,toprule=0.3mm,bottomrule=0.3mm,
  listing options={style=siamlatex}
}

\newtcblisting[use counter=example]{example}[2][]{%
  title={Example~\thetcbcounter: #2},#1}

\newtcbinputlisting[use counter=example]{\examplefile}[3][]{%
  title={Example~\thetcbcounter: #2},listing file={#3},#1}

\DeclareTotalTCBox{\code}{ v O{} }
{ 
  fontupper=\ttfamily\color{black},
  nobeforeafter,
  tcbox raise base,
  colback=codebackground,colframe=white,
  top=0pt,bottom=0pt,left=0mm,right=0mm,
  leftrule=0pt,rightrule=0pt,toprule=0mm,bottomrule=0mm,
  boxsep=0.5mm,
  #2}{#1}

\patchcmd\newpage{\vfil}{}{}{}
\numberwithin{equation}{section}

\newtheorem{theorem}{Theorem}[section]
\newtheorem{lemma}[theorem]{Lemma}
\newtheorem{corollary}[theorem]{Corollary}

\theoremstyle{definition}
\newtheorem{definition}[theorem]{Definition}

\theoremstyle{remark}
\newtheorem{remark}[theorem]{Remark}

\def\rf{\mathrm{f}}

\def\i{\mathrm{i}}

\def\diag{{\mathrm{diag}}}
\def\f{\mathbf{f}}

\def\h{\mathbf{h}}

\def\x{\mathbf{x}}

\def\E{\mathbb{E}}
\def\bS{\mathbb{S}}

\def\N{\mathcal{N}}
\def\M{\mathcal{M}}

\def\0{\mathbf{0}}
\def\1{\mathbf{1}}
\def\V{\mathbb{V}}
\def\E{\mathbb{E}}
\def\F{\mathbb{F}}

\def\I{\mathtt{I}}
\def\B{\mathtt{B}}
\def\beps{{\boldsymbol{\varepsilon}}}
\def\bsigma{{\boldsymbol{\sigma}}}
\def\bfeta{{\boldsymbol{\eta}}}

\def\bfzeta{{\boldsymbol{\zeta}}}

\title[Fundamental Theory and R-linear Convergence of SEM]{Fundamental Theory and R-linear Convergence of Stretch Energy Minimization for Equiareal parameterizations}
\author{Tsung-Ming Huang}
\address{Department of Mathematics, National Taiwan Normal University, Taipei, 116, Taiwan}
\email{min@ntnu.edu.tw}

\author{Wei-Hung Liao}
\address{Department of Applied Mathematics, National Yang Ming Chiao Tung University, Hsinchu, 300, Taiwan}
\email{roger2300245@gmail.com}
\author{Wen-Wei Lin}
\address{Department of Applied Mathematics, National Yang Ming Chiao Tung University, Hsinchu, 300, Taiwan}
\email{wwlin@math.nctu.edu.tw}

\date{\today}

\begin{document}



\maketitle
\begin{abstract}
In this paper, we first extend the finite distortion problem from the bounded domains in $\mathbb{R}^2$ to the closed genus-zero surfaces in $\mathbb{R}^3$ by the stereographic projection. Then we derive a theoretical foundation for spherical equiareal  parameterizations between a simply connected closed surface $\mathcal{M}$ and a unit sphere $\mathbb{S}^2$ via minimizing the total area distortion energy on $\overline{\mathbb{C}}$. Provided we determine the minimizer of the total area distortion energy, the minimizer composed with the initial conformal map determines the equiareal map between the extended planes. Taking the inverse stereographic projection, we can derive the equiareal map between $\mathcal{M}$ and $\mathbb{S}^2$. The total area distortion energy can be rewritten into the sum of Dirichlet energies associated with the southern and northern hemispheres, respectively, and can be decreased by alternatingly solving the corresponding Laplacian equations. Based on this foundational theory, we develop a modified stretch energy minimization for the computation of the spherical equiareal parameterization between $\mathcal{M}$ and $\mathbb{S}^2$. In addition, under some mild conditions, we verify that our proposed algorithm has asymptotically R-linear convergence or forms a quasi-periodic solution. Numerical experiments on various benchmarks validate the assumptions for convergence always hold and indicate the efficiency, reliability and robustness of the developed modified stretch energy minimization.
\end{abstract}


\section{Introduction}
An equiareal (area-preserving) parameterization between a closed genus-zero surface $\mathcal{M}\subset\mathbb{R}^3$ and a unit sphere $\mathbb{S}^2$ has been widely used in various tasks such as surface registration, remeshing, texture mapping and optimal transportation in digital geometry processings. An equiareal parameterization between $\mathcal{M}$ and $\mathbb{S}^2$ can be regarded as to determine an one-to-one correspondence of sampling density on $\mathbb{S}^2$ that minimizes some stretch energy. The total area distortion of the sampling density on $\mathbb{S}^2$ is caused by the stretching of parameterizations. In order to control the density of sampling points on $\mathbb{S}^2$ that minimizes the total area distortion, in this paper, we consider the equiareal parameterization which is one of the best options and refers to as a bijective area-preserving map from $\mathcal{M}$ to $\mathbb{S}^2$.

In fact, there are spherical conformal (angle-preserving) \cite{LHLY} and quasi-conformal parameterizations \cite{ChLL15} between $\mathcal{M}$ and $\mathbb{S}^2$. The uniqueness and existence of the spherical conformal map are claimed by Poincare--Klein--Koebe uniformization Theorem \cite{koebe1907, poincare1908}. The existence of a quasi-conformal mapping is proven in \cite[Chap. 5]{Ahlf66} and the convergence of the iterative algorithm for quasi-conformal mapping is proven in \cite{lui2015}. The conformal parameterization usually minimizes the Dirichlet energy while preserving the minimal angle distortion, i.e., keeping the local shape on the surface as possible up to a similar shape transformation. The quasi-conformal parameterization is used to relax the Beltrami coefficient from zero for conformality to a positive and less than one number. Then, the local shape, e.g., a local disc on the surface can be relaxed to a local ellipse. Conformal or quasi-conformal parameterization can well-preserve the local shapes on the surface. However, in order to achieve conformality or quasi-conformality, it usually maps the area with larger changes on the surface to a smaller area on the sphere, while the area with little changes on the surface occupies a larger area on the sphere. Different from the conformal and quasi-conformal parameterizations with uneven distribution on the sphere, the equiareal parameterization is to consider minimizing the total area distortion, while preserving the local areas on $\mathcal{M}$ and $\mathbb{S}^2$ as possible. The mapping points on the sphere with equiareal parameterization will have a more uniform density distribution.

Recently, efficient algorithms, namely, stretch-minimizing methods \cite{su2016area, yoshizawa2004fast}, Lie advection method \cite{zou2011authalic}, optimal mass transportation methods \cite{TurboSquid, floater2005surface, horn2012matrix, zhao2013area} and stretch energy minimizations (SEM) \cite{YuLL19, yueh2019novel} for the computation of equiareal parameterization have been well developed and utilized in applications. In practice, on effectiveness and accuracy, the SEM algorithm is highly improved compared to the other state-of-the-art algorithms. (see \cite{YuLL19} for details). However, the SEM algorithm still lacks the mathematical foundation for convergence.

In this paper, we aim to provide the fundamental theory for the optimization of the stretching energy. Based on this foundation, we propose the modified SEM for finding the spherical area-preserving parameterization.

The main contributions of this paper are threefold.
\begin{itemize}
\item By considering the stereographic projections and the inversion transformations, we extend the finite distortion problem on bounded domains in $\mathbb{R}^2$ to closed genus-zero surfaces in $\mathbb{R}^3$. We  show that the total area distortion (stretch energy) of $f$ from $\mathcal{M}$ to $\mathbb{S}^2$ is equal to the Dirichlet energy of $g$ on $\overline{\mathbb{C}}$. In particular, we choose a conformal mapping to deform a closed genus-zero surface into a sphere, provided that we derive the minimizer $g$, the mapping $g\circ h$ is the area-preserving map, i.e., $\det(\nabla(g\circ h))=1$.
\item Based on the theoretical foundation above, we propose an efficient and reliable SEM algorithm for the computation of spherical equiareal (area-preserving) parameterization between $\mathcal{M}$ and $\mathbb{S}^2$. Furthermore, we show the SEM algorithm converges asymptotically and R-linearly under some mild conditions or forms a quasi-periodic solution.
\item Numerical experiments on various benchmarks confirm that the assumption for R-linear convergence hold, and the associated means and standard deviation (SD) of local area ratio distortions indicate the accuracy and effectiveness of the proposed SEM algorithm.
\end{itemize}

This paper is organized as follows. In Section~\ref{sec:2}, we introduce the theoretical foundation for the computation of spherical equiareal maps via minimizing the total area distortion energy on the image of $\overline{\mathbb{C}}=\mathbb{R}^2\cup\{\infty\}$. In Section~\ref{sec:3}, we propose a modified SEM for the computation of spherical equiareal parameterizations of simply connected closed triangular meshes and prove the asymptotically R-linear convergence of the SEM algorithm. Numerical experiments on the asymptotically R-linear convergence and the area-preserving property of the SEM are demonstrated in Section \ref{sec:4}. Finally, a concluding remark is given in Section \ref{sec:5}.

In this paper, we use the following notation:
\begin{itemize}[label={$\bullet$}]
\item Bold letters, for instance $\mathbf{f}$, denote vectors in $\mathbb{R}^n$.
\item Capital letters, for instance $A$, denote matrices.
\item Typewriter-style letters, for instance $\mathtt{I}$ and $\mathtt{B}$, denote ordered sets of indices. 
\item $\mathbf{f}_i$ denotes the $i$th entry of the vector $\mathbf{f}$.
\item $\mathbf{f}_\mathtt{I}$ denotes the subvector of $\mathbf{f}$ composed of $\mathbf{f}_i$ for $i\in\mathtt{I}$.
\item $A_{ij}$ denotes the $(i,j)$th entry of matrix $A$.
\item $\lvert A\rvert:=[\lvert A_{ij}\rvert]$ with $\lvert A_{ij}\rvert$ being the modulus of $A_{ij}$.
\item $\bar{A}:=[\bar{A}_{ij}]$ with $\bar{A}_{ij} = \mbox{conj}(A_{ij})$.
\item $A_{\mathtt{I},\mathtt{J}}$ denotes the submatrix of $A$ composed of $A_{ij}$ for $i\in\mathtt{I}$ and $j\in\mathtt{J}$.
\item $\mathbb{S}^2:=\{ \mathbf{x}\in\mathbb{R}^{3} \mid \|\mathbf{x}\|=1 \}$ denotes the $2$-sphere in $\mathbb{R}^{3}$.
\item $\left[{v}_0, \ldots, {v}_m\right]$ denotes the $m$-simplex of the vertices ${v}_0, \ldots, {v}_m$. 
\item $\mathrm{i}$ denotes the imaginary unit $\sqrt{-1}$. 
\end{itemize}

\section{Foundation of Stretch Energy Minimization for Equiareal Map} \label{sec:2}
\subsection{Preliminary}
Let $\Omega$ be an open connected subset of $\mathbb{R}^n$ and     
$f=\left(f^1, \cdots, f^n\right)$ be a map from $\Omega$ to $\mathbb{R}^n$ with  $f^i\in W^{1, n}(\Omega)$ for $i=1, \ldots, n$.

\begin{definition}\cite{astala2005extremal}
\label{def:finite1}
Suppose $f\in W^{1, n}_{loc}\left(\Omega, \mathbb{R}^n\right)$ with nonnegative Jacobian 
$J(\mathbf{x}, f)$. $f$ is said to have the finite distortion $K$, if there exists a function $K(\mathbf{x})\in [1, \infty)$ defined a.e. in $\Omega$ such that
\begin{equation*}
|\nabla f(\mathbf{x})|_{\delta}^2=K(\mathbf{x})J(\mathbf{x}, f), \quad \mbox{a.e. $\mathbf{x}\in\Omega$},
\end{equation*}
where $|\nabla f(\mathbf{x})|_{\delta}^2:=\mathrm{Tr}\left({\nabla f(\mathbf{x})}^{\top}\nabla  f(\mathbf{x})\right)$ is the Euclidean norm of the differential $\nabla f:\mathbb{R}^n\to\mathbb{R}^n$ and $J(\mathbf{x}, f)=\det \nabla f$. Here, ``$\mathrm{Tr}$'' is ``\textup{Trace}'' for short.
\end{definition}
Note that the weak assumption $ f\in W^{1, n}_{loc}\left(\Omega, \mathbb{R}^n\right)$ does not ensure that $ f$ is a continuous mapping, but $ f$ is continuous if $ f$ has the finite distortion \cite{vodop1976}.

At first, a manifold $\mathcal{M}^n$ of dimensional $n$ is a locally Euclidean topological space, that is, there is a neighborhood around any point in $\mathcal{M}^n$ which is topologically the same as the open set in $\mathbb{R}^n$. If we consider some local properties of $ f$ that are related to the finite distortion, then the above arguments can be directly applied on any neighborhood of $\mathcal{M}^n$. In other words, the metrics in either the domain or the image can be chosen to be the standard Euclidean metric which is not associated with $ f$. Once we want to analyze some global properties of $ f$ with the finite distortion, it is inevitable to consider globally defined metrics in either the domain or the image.

In general, let $(\mathcal{M}, \gamma)$ and $(\mathcal{N}, \rho)$ be Riemannian manifolds of dimension two with $\gamma:=[\gamma_{\alpha\beta}]_{\alpha, \beta= 1, 2}$ and $\rho:=[\rho_{ij}]_{i, j= 1, 2}$ being the associated metrics, respectively. Consider a Sobolev mapping $ f:\mathcal{M}\to\mathcal{N}$ in $W^{1, 2}(\mathcal{M}, \mathcal{N})$. We extend Definition~\ref{def:finite1} on manifolds.
\begin{definition}
\label{def:finite2}
Suppose $ f\in W^{1, 2}\left(\mathcal{M}, \mathcal{N}\right)$ with nonnegative Jacobian $J(\mathbf{x},  f)_{\rho}$.  $ f$ is said to have the finite distortion $K$, if there exists a function $K(\mathbf{x})\in [1, \infty)$ defined a.e. in $\mathcal{M}$ such that
\begin{equation*}
|\nabla f(\mathbf{x})|_{\rho}^2=K(\mathbf{x})J(\mathbf{x},  f)_{\rho}, \quad \mbox{a.e. $\mathbf{x}\in\mathcal{M}$},  
\end{equation*}
where the metric norm of the differential $\nabla f:\mathbb{R}^2\to\mathbb{R}^2$ is defined as
\begin{equation}
\label{grad}
|\nabla f(\mathbf{x})|^2_{\rho}:=\gamma^{\alpha\beta}{\rho}_{ij}( f(\mathbf{x}))\frac{\partial f^i}{\partial x^{\alpha}}\frac{\partial f^j}{\partial x^{\beta}},
\end{equation}
and the nonnegative Jacobian $J(\mathbf{x},  f)_{\rho}$ is defined as
\begin{align}
\label{jaco}
J(\mathbf{x},  f)_{\rho}:=\mbox{det}_{\rho}\nabla f=\left(\det\left({\nabla f(\mathbf{x})}^{\top}\cdot_{\rho}\nabla f(\mathbf{x})\right)\right)^{1/2}
\end{align}
in which $\cdot_{\rho}$ is the inner product induced by the associated metric $\rho$. Note that in \eqref{grad} we express the summation by using the conventional Einstein notations and $[\gamma^{\alpha\beta}]  :=[\gamma_{\alpha\beta}]^{-1}$.
\end{definition}

The following theorem states the invertibility of the Sobolev functions.
\begin{theorem}\cite{fonseca1995local}
\label{inversobolev}
Let $\Omega\subset\mathbb{R}^n$ be a bounded, open set and let $ f\in W^{1, n}(\Omega, \mathbb{R}^n)$ be a function such that $\det\nabla  f(\mathbf{x})>0$ a.e. $\mathbf{x}\in\Omega$. Then for almost every $\mathbf{x}_0\in\Omega$, $ f$ is locally almost invertible in a neighborhood of $\mathbf{x}_0$, in the sense that there exists $r=r(\mathbf{x}_0)>0$, an open set $D=D(\mathbf{x}_0)$, and a function $ g:B(\mathbf{y}_0, r)\to D$ with $\mathbf{y}_0= f(\mathbf{x}_0)$ such that $ g\in W^{1, 1}(B(\mathbf{y}_0, r), \mathbb{R}^n)$, $ g\circ  f(\mathbf{x})=\mathbf{x}$, $ f\circ  g(\mathbf{y})=\mathbf{y}$ and $\nabla  g(\mathbf{y})=(\nabla  f)^{-1}( g(\mathbf{y}))$ for a.e. $\mathbf{x}\in D$ and $\mathbf{y}\in B(\mathbf{y}_0, r).$ If, in addition, ${\lvert\frac{\textup{adj}\left(\nabla  f\right)}{\det\nabla f}\rvert}^s\det\nabla f\in L^1(\Omega)$ for some $1\leq s<\infty$, then $ g\in W^{1, s}\left(B(\mathbf{y}_0, r), D\right)$.
\end{theorem}

\subsection{Finite Distortion between Compact Manifolds}
Let $(\mathbb{S}^2, \delta)$ be a unit sphere equipped with the standard Euclidean metric $\delta:=\delta_{ij} dx^i dx^j$, $i, j=1, 2, 3$, where $\delta_{ij}$ is the Kronecker delta, and $(\mathcal{M},  f^*\delta)$ be a closed genus-zero surface in $\mathbb{R}^3$ with the pullback metric $ f^*\delta$, where $ f\in W^{1, 2}(\mathcal{M}, \mathbb{R}^3)$ such that $\det\nabla  f(\mathbf{x})>0$ a.e. $\mathbf{x}\in\mathcal{M}$. 
We now use the stereographic projection to transform the problem in $\mathbb{R}^3$ together with their induced metrics into the problem in $\mathbb{R}^2$.

Let $(\Omega_0, \delta)$ and $(\Omega, \rho)$ be the subsets in $\pi(\mathcal{M})$ and $\pi(\mathbb{S}^2)$ with induced metrics $\delta:=\delta_{ij}dx^idx^j, i, j=1, 2$ and $\rho:=\frac{4}{(u^2+v^2+1)^2}\left(du^2+dv^2\right)$, respectively, where
\begin{align*}
    \pi:=\Pi_{\mathbb{S}^2}:\mathbb{S}^2\longrightarrow\mathbb{R}^2,\quad(x^1, x^2, x^3)\longmapsto\left(\frac{x^1}{1-x^3}, \frac{x^2}{1-x^3}\right)
\end{align*}
denotes the stereographic projection.
Let
\begin{align}\label{eq:h}
    h:=\pi\circ f\circ\pi^{-1}
\end{align}
in $W^{1, 2}(\Omega_0, \Omega)$ be the Sobolev map with $\det\nabla  h(\mathbf{x})>0$ a.e. $\mathbf{x}\in\Omega_0$. According to Definition~\ref{def:finite2}, we compute the finite distortion $K$ of $ h=(h^1, h^2)$ as
\begin{align}
\frac{|\nabla h|_{\rho}^2}{J(\mathbf{x},  h)_{\rho}}&=\frac{\delta^{\alpha\beta}\rho_{ij}( h(\mathbf{x}))\frac{\partial h^i}{\partial x^{\alpha}}\frac{\partial h^j}{\partial x^{\beta}}}{\sqrt{\det(\nabla h^{\top}\cdot_{\rho} \nabla h)}}\notag\\
&=\frac{\text{Tr}\left(
\begin{bmatrix}
\frac{\partial h^1}{\partial x^{1}}&\frac{\partial h^2}{\partial x^{1}}\\
\frac{\partial h^1}{\partial x^{2}}&\frac{\partial h^2}{\partial x^{2}}
\end{bmatrix}
\begin{bmatrix}
\frac{4}{(| h|^2+1)^2}&0\\
0&\frac{4}{(| h|^2+1)^2}
\end{bmatrix}
\begin{bmatrix}
\frac{\partial h^1}{\partial x^{1}}&\frac{\partial h^1}{\partial x^{2}}\\
\frac{\partial h^2}{\partial x^{1}}&\frac{\partial h^2}{\partial x^{2}}
\end{bmatrix}\right)
}{\sqrt{\det\left(\begin{bmatrix}
\frac{\partial h^1}{\partial x^{1}}&\frac{\partial h^2}{\partial x^{1}}\\
\frac{\partial h^1}{\partial x^{2}}&\frac{\partial h^2}{\partial x^{2}}
\end{bmatrix}
\begin{bmatrix}
\frac{4}{(| h|^2+1)^2}&0\\
0&\frac{4}{(| h|^2+1)^2}
\end{bmatrix}
\begin{bmatrix}
\frac{\partial h^1}{\partial x^{1}}&\frac{\partial h^1}{\partial x^{2}}\\
\frac{\partial h^2}{\partial x^{1}}&\frac{\partial h^2}{\partial x^{2}}
\end{bmatrix}\right)}}\notag\\
&=\frac{|\nabla h|_{\delta}^2}{J(\mathbf{x},  h)}.\label{eq:distquo}
\end{align}
We observe that the finite distortion in Definition~\ref{def:finite2} is coincident with the finite distortion in Definition~\ref{def:finite1} if the image of $h$ is endowed with a conformal metric.

\begin{definition}
Let $ h\in W^{1, 2}(\Omega_0, \Omega)$  as in \eqref{eq:h}  with $\det\nabla h(\mathbf{x})>0$ a.e. $\mathbf{x}\in\Omega_0$. The total area distortion (stretch energy) of $ h$ on $\Omega_0$ is defined as 
\begin{align}
    E_{\rho}( h; \Omega_0)=\int_{\Omega_0}\frac{|\nabla h|_{\rho}^2}{J(\mathbf{x},  h)_{\rho}} dx^1\wedge dx^2.\label{totaldis}
\end{align}
\end{definition}

\begin{lemma}
\label{obser}
Let $\Omega_0$ be a bounded domain in $\mathbb{R}^2$ with metric $\delta$ and $ h$ in \eqref{eq:h} be a mapping in $W^{1, 2}(\Omega_0, \mathbb{R}^2)$ with $\det\nabla  h(\mathbf{x})>0$ a.e. $\mathbf{x}\in\Omega_0$. Then $ g:= h^{-1}$ exists in $W^{1, 2}(\Omega, \mathbb{R}^2)$ with  $\Omega:= h(\Omega_0)$ and the total area distortion $E_{\rho}( h; \Omega_0)$ in \eqref{totaldis} is equal to the Dirichlet energy of $ g$ on $\Omega$, i.e.,
\begin{align}
\label{eq:integral}
E_{\rho}( h; \Omega_0)
=\int_{\Omega}|\nabla_{\mathbf{y}} g(\mathbf{y})|_{\delta}^2 dy^1\wedge dy^2\equiv 
E( g; \Omega).
\end{align}
\end{lemma}
\begin{proof}
Since $ h$ is a mapping  in $W^{1, 2}(\Omega_0, \Omega)$ and $\det\nabla h(\mathbf{x})>0$ a.e. $\mathbf{x}\in\Omega_0$, the existence of $ g\in W^{1, 2}(\Omega, \mathbb{R}^2)$ is guaranteed by Theorem~\ref{inversobolev}. Therefore, we only need to check the total area distortion \eqref{eq:integral},
\begin{align*}
E_{\rho}( h; \Omega_0)
&\equiv\int_{\Omega_0}\frac{|\nabla_{\mathbf{x}} h(\mathbf{x})|_{\rho}^{2}}{J(\mathbf{x},  h)_{\rho}}dx^1\wedge dx^2\\
&=\int_{\Omega_0}\frac{|\nabla_{\mathbf{x}} h(\mathbf{x})|_{\delta}^{2}}{J(\mathbf{x},  h)}dx^1\wedge dx^2\\
&=\int_{\Omega_0}\frac{|\textup{adj}(\nabla_{\mathbf{x}} h(\mathbf{x}))|_{\delta}^2}{J(\mathbf{x},  h)}dx^1\wedge dx^2\\
&=\int_{\Omega_0}|\nabla_{\mathbf{x}} h^{-1}|_{\delta}^2 J(\mathbf{x},  h)dx^1\wedge dx^2\quad (\mathbf{y}= h(\mathbf{x}))\\
&=\int_{\Omega}|\nabla_{\mathbf{y}} g(\mathbf{y})|_{\delta}^2dy^1\wedge dy^2.
\end{align*}
Here, the first three equalities are derived by \eqref{eq:distquo} and the equality $M \textup{adj}(M)= \textup{adj}(M) M=\det(M) I$, where $M$ is a matrix.
\end{proof}

Based on Lemma~\ref{obser}, the original problem is illustrated as follows.
Given a mapping $ h$ in $W^{1, 2}(\Omega_0, \Omega)$ and $\det\nabla h(\mathbf{x})>0$ a.e. $\mathbf{x}\in\Omega_0$, the minimization of the total area distortion $E_{\rho}( h; \Omega_0)$ in \eqref{eq:integral} is equivalent to minimizing the Dirichlet energy $E( g; \Omega)$, which becomes to solve $\Delta_\mathbf{y}  g=0$.

\begin{remark}
\begin{enumerate}
\item[(i)] \label{1}Area distortion.
We first apply the stereographic projection $\pi$ to flatten a closed genus-zero surface onto a plane. This step causes the area distortion rate, say $\mu$. Next, we use the spherical CEM algorithm in \cite{LHLY} to obtain a conformal mapping $ h\rvert_{\Omega_0}:\Omega_0\to\Omega$ , $( h\rvert_{\Omega_0^{\complement}}:\Omega_0^{\complement}\to\Omega^{\complement})$ which produces another area distortion rate, say $\nu$. Then, we consider the minimization of the total area distortion problem and determine the conformal mapping $ g\rvert_{\Omega}:\Omega\to g(\Omega)$, $( g\rvert_{\Omega^{\complement}}:\Omega^{\complement}\to g(\Omega^{\complement}))$ which modifying the area distortion rate by $1/\nu$. Finally, we apply the inverse stereographic projection $\pi^{-1}$ to map $\overline{\mathbb{R}}^2$ to a unit sphere $\mathbb{S}^2$, which amends the area distortion in the scale of $1/\mu$. Hence, the total area distortion rate is 1 through the mapping $\pi^{-1}\circ g\circ  h\circ\pi$. Namely, we obtain an area-preserving map between $\mathcal{M}$ and $\mathbb{S}^2$.
\item[(ii)] Angle distortion.
The angle distortion happens in the last step that the inverse stereographic projection maps $\overline{\mathbb{R}}^2$ back to $\mathbb{S}^2$, since the Dirichlet energy $E( g; \Omega)$ and $E( g; \Omega^{\complement})$ are not the Dirichlet energy determined by $\pi$. In other words, $E( g; \Omega)$ and $E( g; \Omega^{\complement})$ should be given by
\begin{align*}
\int_{\Omega}|\nabla_{\mathbf{y}} g(\mathbf{y})|_{\rho}^2dy^1\wedge dy^2\equiv
\int_{\Omega}\frac{4}{\left(1+|g(\mathbf{y})|\right)^2}|\nabla_{\mathbf{y}} g(\mathbf{y})|_{\delta}^2dy^1\wedge dy^2,
\end{align*}
and 
\begin{align*}
\int_{\Omega^{\complement}}|\nabla_{\mathbf{y}} g(\mathbf{y})|_{\rho}^2dy^1\wedge dy^2\equiv    
\int_{\Omega^{\complement}}\frac{4}{\left(1+|g(\mathbf{y})|\right)^2}|\nabla_{\mathbf{y}} g(\mathbf{y})|_{\delta}^2dy^1\wedge dy^2.
\end{align*}
Namely, the image of $g$ should be accompanied by the metric $\rho$ rather than $\delta$.
\end{enumerate}
\end{remark}
The conclusion of two distortions is meaningful because the maps that preserve angles and areas are rigid transformations such as rotations and reflections. However, the rigid transformation cannot deform a closed non-sphere into a sphere.

\subsection{Finite Distortion between Closed Manifolds}
Let $ f:\mathcal{M}\to\mathbb{S}^2$ be a map in $W^{1, 2}(\mathcal{M}, \mathbb{S}^2)$ with $\det\nabla f(\mathbf{x})>0$ a.e. $\mathbf{x}\in\mathcal{M}$. 
We transform this problem onto the extended plane by applying the stereographic projection $\pi$. Thus, we alternatively consider the map $h$ of \eqref{eq:h} in $W^{1, 2}(\pi(\mathcal{M}), \overline{\mathbb{R}}^2)$  with $\det\nabla h(\mathbf{x})>0$, and then we divide $ h$ into $ h_{S}$ and $ h_{N}$,
\begin{align}
\label{cfmap}
 h(\mathbf{x})=
\begin{cases}
 h_{S}(\mathbf{x})\in \mathbb{D},& \mbox{if $\mathbf{x}\in\Omega_0$,}\\
 h_{N}(\mathbf{x})\in \overline{\mathbb{R}}^2\backslash\mathbb{D},& \mbox{if $\mathbf{x}\in\overline{\mathbb{R}}^2\backslash\Omega_0$,}
\end{cases}
\end{align}
where $\mathbb{D}$ is a unit disk on $\mathbb{R}^2$ and $\Omega_0$ is some open subset in $\mathbb{R}^2$.
\begin{definition}
Let $ f:\mathcal{M}\to\mathbb{S}^2$ be a map in $W^{1, 2}(\mathcal{M}, \mathbb{S}^2)$ with $\det\nabla f(\mathbf{x})>0$ a.e. $\mathbf{x}\in\mathcal{M}$ and $ h$ be given in \eqref{cfmap}. The total area distortion (stretch energy) of $ h$ on $\overline{\mathbb{R}}^2$ is defined as 
\begin{equation}\label{Totaldis}
    E_{\rho}( h; \overline{\mathbb{R}}^2)=\int_{\overline{\mathbb{R}}^2}\frac{|\nabla h(\mathbf{x})|_{\rho}^{2}}{J(\mathbf{x},  h)_{\rho}}dx^1\wedge dx^2.
\end{equation}
\end{definition}
Because the representation of the northern and southern hemispheres depends on the chosen coordinate plane, we can re-express the function $ h$ as a mapping between bounded regions 
\begin{align}
\label{recfmap}
 h(\mathbf{x})=
\begin{cases}
h_{1}(\mathbf{x})\in \mathbb{D},& \mbox{if $\mathbf{x}\in\Omega_0$,}\\
h_{2}(\hat{\mathbf{x}})\in \mathbb{D},& \mbox{if $\hat{\mathbf{x}}:=\frac{\mathbf{x}}{|\mathbf{x}|^2}\in\hat{\Omega}_0$.}
\end{cases}
\end{align}
where $\hat{\Omega}_0:=\textup{Inv}(\overline{\mathbb{R}}^2\backslash\Omega_0)$, $h_{2}:=\textup{Inv}( h_{N})$, and $h_{1}:= h_{S}$.
\begin{theorem}
\label{invconformal}
Let $(\pi(\mathcal{M}), \delta)$ and $(\overline{\mathbb{R}}^2, \rho)$ be two Riemannian manifolds with metrics $\delta:=\delta_{ij}dx^idx^j, i, j=1, 2$ and $\rho:=\frac{4}{(u^2+v^2+1)^2}\left(du^2+dv^2\right)$. Let $ h$ in \eqref{eq:h} be a mapping in $W^{1, 2}(\pi(\mathcal{M}), \overline{\mathbb{R}}^2)$ with $\det\nabla h(\mathbf{x})>0$. Then $ g:= h^{-1}$ exists in $W^{1, 2}(\overline{\mathbb{R}}^2, \overline{\mathbb{R}}^2)$ and is equal to
\begin{align}
\label{invh}
 g(\mathbf{y})
&=\mathbbm{1}_{\mathbb{D}}(\mathbf{y})h_{1}^{-1}(\mathbf{y})+\mathbbm{1}_{\mathbb{D}}(\frac{\mathbf{y}}{|\mathbf{y}|^2})h_{2}^{-1}(\frac{\mathbf{y}}{|\mathbf{y}|^2})\\
&:=g_{1}(\mathbf{y})+g_{2}(\hat{\mathbf{y}})\nonumber, 
\end{align}
where $\mathbbm{1}_{\mathbb{D}}(\cdot)$ is the characteristic function on $\mathbb{D}$. Furthermore, the total area distortion \eqref{Totaldis} of $ h$ is equal to the Dirichlet energy $E( g)$ of $ g$,
\begin{align}
\int_{\overline{\mathbb{R}}^2}\frac{|\nabla h(\mathbf{x})|_{\rho}^{2}}{J(\mathbf{x},  h)_{\rho}}dx^1\wedge dx^2
&=\int_{\mathbb{D}}|\nabla_{\mathbf{y}}g_{1}(\mathbf{y})|_{\delta}^2dy^1\wedge dy^2+\int_{\mathbb{D}}|\nabla_{\hat{\mathbf{y}}}g_{2}(\hat{\mathbf{y}})|_{\delta}^2 d\hat{y}^2\wedge d\hat{y}^1\equiv E( g),\label{finite2}
\end{align}
where $\mathbf{y}=h_{1}(\mathbf{x})$ and $\hat{\mathbf{y}}=h_{2}(\hat{\mathbf{x}})$ are defined by \eqref{recfmap}.
\end{theorem}
\begin{proof}
First, we separate the total area distortion \eqref{Totaldis} into two terms,
\begin{align}
\label{Totaldis1}
\int_{\overline{\mathbb{R}}^2}\frac{|\nabla h(\mathbf{x})|_{\rho}^{2}}{J(\mathbf{x},  h)_{\rho}}dx^1\wedge dx^2
&=\int_{\Omega_0}\frac{|\nabla h_S(\mathbf{x})|_\rho^{2}}{J(\mathbf{x},  h_S)_\rho}dx^1\wedge dx^2+\int_{\overline{\mathbb{R}}^2\backslash\Omega_0}\frac{|\nabla h_N(\mathbf{x})|_{\rho}^{2}}{J(\mathbf{x},  h_N)_{\rho}}dx^1\wedge dx^2.
\end{align}
(i) Based on Lemma~\ref{obser}, the total area distortion of $ h_{S}$ in $\Omega_0$ is well-defined and then the first term in \eqref{Totaldis1} is equivalent to the Dirichlet energy of $g_{1}:=h_{1}^{-1}$ in $\mathbb{D}$, i.e., 
\begin{align}
\int_{\Omega_0}\frac{|\nabla_{\mathbf{x}} h_{S}(\mathbf{x})|_{\rho}^{2}}{J(\mathbf{x},  h_{S})_{\rho}}dx^1\wedge dx^2
&=\int_{\Omega_0}\frac{|\nabla_{\mathbf{x}}h_{1}(\mathbf{x})|_{\rho}^{2}}{J(\mathbf{x}, h_{1})_{\rho}}dx^1\wedge dx^2\quad (\text{by \eqref{cfmap} and \eqref{recfmap}})\nonumber\\
&=\int_{\mathbb{D}}|\nabla_{\mathbf{y}}g_{1}(\mathbf{y})|_{\delta}^2dy^1\wedge dy^2.\quad (\mbox{by Lemma~\ref{obser}}).\label{totaldis:south}
\end{align}
Here, $g_{1}:=h_{1}^{-1}$ defines the first part in \eqref{invh}.

\noindent (ii) Since the northern and southern hemispheres are equivalent to each other, the second term in \eqref{Totaldis1}  can be modified to the total area distortion problem between bounded domains by the inversion transformation. Applying the local invertibility of Theorem~\ref{inversobolev}, \eqref{grad} and \eqref{jaco}, we have 
\begin{align}
\frac{|\nabla_{\mathbf{x}} h_{N}(\mathbf{x})|_{\rho}^2}{J(\mathbf{x},  h_{N})_\rho}
&=\frac{\delta^{\alpha\beta}\rho_{ij}( h_{N}(\mathbf{x}))\frac{\partial h_{N}^i}{\partial x^{\alpha}}\frac{\partial h_{N}^j}{\partial x^{\beta}}}{\sqrt{\det(\nabla_{\mathbf{x}} h_{N}^{\top}\cdot_{\rho} \nabla_{\mathbf{x}} h_{N})}}\notag\\
&=\frac{\mathrm{Tr}\left(
\begin{bmatrix}
1&0\\
0&1
\end{bmatrix}
\begin{bmatrix}
\frac{\partial h_{N}^1}{\partial x^{1}}&\frac{\partial h_{N}^2}{\partial x^{1}}\\
\frac{\partial h_{N}^1}{\partial x^{2}}&\frac{\partial h_{N}^2}{\partial x^{2}}
\end{bmatrix}
\begin{bmatrix}
\frac{4}{(| h_{N}|^2+1)^2}&0\\
0&\frac{4}{(| h_{N}|^2+1)^2}
\end{bmatrix}
\begin{bmatrix}
\frac{\partial h_{N}^1}{\partial x^{1}}&\frac{\partial h_{N}^1}{\partial x^{2}}\\
\frac{\partial h_{N}^2}{\partial x^{1}}&\frac{\partial h_{N}^2}{\partial x^{2}}
\end{bmatrix}\right)
}{\sqrt{\det\left(
\begin{bmatrix}
\frac{\partial h_{N}^1}{\partial x^{1}}&\frac{\partial h_{N}^2}{\partial x^{1}}\\
\frac{\partial h_{N}^1}{\partial x^{2}}&\frac{\partial h_{N}^2}{\partial x^{2}}
\end{bmatrix}
\begin{bmatrix}
\frac{4}{(| h_{N}|^2+1)^2}&0\\
0&\frac{4}{(| h_{N}|^2+1)^2}
\end{bmatrix}
\begin{bmatrix}
\frac{\partial h_{N}^1}{\partial x^{1}}&\frac{\partial h_{N}^1}{\partial x^{2}}\\
\frac{\partial h_{N}^2}{\partial x^{1}}&\frac{\partial h_{N}^2}{\partial x^{2}}
\end{bmatrix}\right)}}\notag\\
&\equiv \frac{\mathcal{P}}{\mathcal{Q}}.\label{quotient:north}
\end{align}
With $h_{2}:=\frac{1}{\overline{ h_{N}}}$ and the identity 
\begin{align*}
\langle \nabla_{\mathbf{x}}\frac{1}{\overline{ h_{N}}}, \nabla_{\mathbf{x}}\frac{1}{\overline{ h_{N}}}\rangle=
&\sum_{i=1}^2\left( \left(\frac{-\partial h_{N}^1/\partial x^i}{| h_{N}|^2}\right)^2+\left(\frac{-\partial h_{N}^2/\partial x^i}{| h_{N}|^2}\right)^2\right)\nonumber\\
=&\frac{1}{| h_{N}|^4}\langle \nabla_{\mathbf{x}} h_{N}, \nabla_{\mathbf{x}} h_{N}\rangle,
\end{align*}
we rewrite $\mathcal{P}$ in \eqref{quotient:north} into the form of $h_{2}$ as follows,
\begin{align}
\label{finite:north1}
\mathcal{P}=\frac{4| h_{N}|^4}{(1+| h_{N}|^2)^2}\langle \nabla_{\mathbf{x}}\frac{1}{\overline{ h_{N}}}, \nabla_{\mathbf{x}}\frac{1}{\overline{ h_{N}}}\rangle=\frac{4}{(1+|h_{2}|^2)^2}\langle \nabla_{\mathbf{x}}h_{2}, \nabla_{\mathbf{x}}h_{2}\rangle.
\end{align}
On the other hand, since
\begin{align*}
    A &\equiv \begin{bmatrix}
\frac{\partial h_{N}^1}{\partial x^{1}}&\frac{\partial h_{N}^1}{\partial x^{2}}\\
\frac{\partial h_{N}^2}{\partial x^{1}}&\frac{\partial h_{N}^2}{\partial x^{2}}\\
\end{bmatrix} 
    = {\large \begin{bmatrix}
\frac{\frac{\partial h^1_2}{\partial x^1}|h_{2}|^2-h_{2}^1(2h_{2}^1\frac{\partial h^1_2}{\partial x^1}+2h_{2}^2\frac{\partial h^2_2}{\partial x^1})}{|h_{2}|^4}&\frac{\frac{\partial h^1_2}{\partial x^2}|h_{2}|^2-h_{2}^1(2h_{2}^1\frac{\partial h^1_2}{\partial x^2}+2h_{2}^2\frac{\partial h^2_2}{\partial x^2})}{|h_{2}|^4}\\
\frac{\frac{\partial h^2_2}{\partial x^1}|h_{2}|^2-h_{2}^2(2h_{2}^1\frac{\partial h^1_2}{\partial x^1}+2h_{2}^2\frac{\partial h^2_2}{\partial x^1})}{|h_{2}|^4}&\frac{\frac{\partial h^2_2}{\partial x^2}|h_{2}|^2-h_{2}^2(2h_{2}^1\frac{\partial h^1_2}{\partial x^2}+2h_{2}^2\frac{\partial h^2_2}{\partial x^2})}{|h_{2}|^4}
\end{bmatrix}}, \nonumber \\
D_\rho &\equiv \begin{bmatrix}
\frac{4}{(| h_{N}|^2+1)^2}&0\\
0&\frac{4}{(| h_{N}|^2+1)^2}\\
\end{bmatrix} = \begin{bmatrix} 
     \frac{4|h_{2}|^4}{(1+|h_{2}|^2)^2} & 0 \\ 0 & \frac{4|h_{2}|^4}{(1+|h_{2}|^2)^2}
\end{bmatrix},
\end{align*}
it holds that
\begin{align}
\label{finite:north4}
\det(A)=\frac{|h_{2}|^4(\frac{\partial h^2_2}{\partial x^1}\frac{\partial h^1_2}{\partial x^2}-\frac{\partial h^1_2}{\partial x^1}\frac{\partial h^2_2}{\partial x^2})}{|h_{2}|^8},\quad \det(D_rho)=\frac{16|h_{2}|^8}{(1+|h_{2}|^2)^4}.
\end{align}
Because the inversion transformation is orientation-reversing, from \eqref{finite:north4} we represent $\mathcal{Q}$ in \eqref{quotient:north} by $h_{2}$ as,
\begin{align}
\label{finite:north7}
\mathcal{Q}=-\sqrt{ \frac{16|h_{2}|^8(\frac{\partial h^2_2}{\partial x^1}\frac{\partial h^1_2}{\partial x^2}-\frac{\partial h^1_2}{\partial x^1}\frac{\partial h^2_2}{\partial x^2})^2}{(1+|h_{2}|^2)^4|h_{2}|^8}}=-\frac{4}{(1+|h_{2}|^2)^2}J(x, h_{2}).
\end{align}
Combining \eqref{finite:north1} with \eqref{finite:north7}, the second term in \eqref{Totaldis1} can be rewritten as the total area distortion of $h_{2}$,
\begin{align}
\int_{\overline{\mathbb{R}}^2\backslash\Omega_0}\frac{|\nabla_{\mathbf{x}} h_{N}(\mathbf{x})|_{\rho}^{2}}{J(\mathbf{x},  h_{N})_{\rho}}dx^1\wedge dx^2
&=\int_{\overline{\mathbb{R}}^2\backslash\Omega_0}\frac{|\nabla_{\mathbf{x}}h_{2}(\mathbf{x})|_{\rho}^{2}}{-J(\mathbf{x}, h_{2})_{\rho}}dx^1\wedge dx^2\nonumber\\
&=\int_{\overline{\mathbb{R}}^2\backslash\Omega_0}\frac{|\nabla_{\mathbf{x}}h_{2}(\mathbf{x})|_{\delta}^{2}}{-J(\mathbf{x}, h_{2})}dx^1\wedge dx^2.\quad (\mbox{by \eqref{eq:distquo}})\label{northNinv}
\end{align}
Next, we consider the coordinate transformation $\mathcal{T}$,
\begin{align}
\label{coord}
\mathcal{T}(\overline{\mathbb{R}}^2\backslash\Omega_0)=\hat{\Omega}_0,\quad \mathcal{T}(\mathbf{x})=(\frac{x^1}{|\mathbf{x}|^2}, \frac{x^2}{|\mathbf{x}|^2}):=(\hat{x}^1, \hat{x}^2)=\hat{\mathbf{x}}.
\end{align}
For clarity, we change the upper index to the lower index. In other words, $x^i$ is the $i$th component of the vector $\mathbf{x}$, but we will denote $x_i^k:=(x^i)^k$ be the $k$ power of the $i$th component. Then  the 1-forms, 2-form, and the partial derivatives of $h_{2}$ in \eqref{northNinv} under the transformation \eqref{coord} can be expressed as follows,
\begin{subequations}
\begin{align}
\begin{bmatrix}
d\hat{x}^1\\
d\hat{x}^2
\end{bmatrix}
&=
\begin{bmatrix}
\frac{x_2^2-x_1^2}{(x_1^2+x_2^2)^2}&\frac{-2x_1x_2}{(x_1^2+x_2^2)^2}\\
\frac{-2x_1x_2}{(x_1^2+x_2^2)^2}&\frac{x_1^2-x_2^2}{(x_1^2+x_2^2)^2}
\end{bmatrix}
\begin{bmatrix}
dx^1\\
dx^2
\end{bmatrix},\label{1form} \\ 
d\hat{x}^1\wedge d\hat{x}^2 &=-\frac{1}{(x_1^2+x_2^2)^2}dx^1\wedge dx^2,\label{2form} \\ 
\nabla_{\mathbf{x}}h_{2}(\mathbf{x})
&=\begin{bmatrix}
\frac{\partial h_{2}^1}{\partial x^1}&\frac{\partial h_{2}^1}{\partial x^2}\\
\frac{\partial h_{2}^2}{\partial x^1}&\frac{\partial h_{2}^2}{\partial x^2}
\end{bmatrix}
=\begin{bmatrix}
\frac{\partial h_{2}^1}{\partial \hat{x}^1}&\frac{\partial h_{2}^1}{\partial \hat{x}^2}\\
\frac{\partial h_{2}^2}{\partial \hat{x}^1}&\frac{\partial h_{2}^2}{\partial \hat{x}^2}
\end{bmatrix}
\begin{bmatrix}
\frac{\partial \hat{x}^1}{\partial x^1}&\frac{\partial \hat{x}^1}{\partial x^2}\\
\frac{\partial \hat{x}^2}{\partial x^1}&\frac{\partial \hat{x}^2}{\partial x^2}
\end{bmatrix}
=\nabla_{\hat{\mathbf{x}}}h_{2}\nabla_{\mathbf{x}}\hat{\mathbf{x}}.\label{deriva}
\end{align}
\end{subequations}
According to Definition~\ref{def:finite1}, \eqref{deriva} and \eqref{coord}, it holds that
\begin{subequations} \label{eq:J_gradf}
\begin{align}
J(\mathbf{x}, h_{2})	&=\det \nabla_{\hat{\mathbf{x}}}h_{2}\ \det \nabla_{\mathbf{x}}\hat{\mathbf{x}}=-\frac{1}{(x_1^2+x_2^2)^2}J(\hat{\mathbf{x}}, h_{2}), \label{Jaco} \\
|\nabla_{\mathbf{x}}h_{2}(\mathbf{x})|_{\delta}^2 &= \mathrm{Tr}\left((\nabla_{\hat{\mathbf{x}}}h_{2}\nabla_{\mathbf{x}}\hat{\mathbf{x}})^{\top}(\nabla_{\hat{\mathbf{x}}}h_{2}\nabla_{\mathbf{x}}\hat{\mathbf{x}})\right)=\frac{1}{(x_1^2+x_2^2)^2}|\nabla_{\hat{\mathbf{x}}}h_{2}(\hat{\mathbf{x}})|_{\delta}^2\label{Diri}
\end{align} 
\end{subequations}
Using\eqref{coord}, \eqref{2form} and \eqref{eq:J_gradf}, the total area distortion in \eqref{northNinv}  can be rewritten as
\begin{align}
\label{Totaldis:north}
\int_{\overline{\mathbb{R}}^2\backslash\Omega_0}\frac{|\nabla_{\mathbf{x}}h_{2}(\mathbf{x})|_{\delta}^{2}}{-J(\mathbf{x}, h_{2})}dx^1\wedge dx^2
&=\int_{\hat{\Omega}_0}\frac{|\nabla_{\hat{\mathbf{x}}}h_{2}(\hat{\mathbf{x}})|_{\delta}^2}{J(\hat{\mathbf{x}}, h_{2})}d\hat{x}^2\wedge d\hat{x}^1.
\end{align}
Hence, the total area distortion of $ h$ between unbounded domains $\Omega_0^{\complement}$ and $\mathbb{D}^{\complement}$ is turned into the total area distortion of $\hat{ f}$ between bounded domains $\hat{\Omega}_0$ and $\mathbb{D}$. 
We use Lemma~\ref{obser} to derive $g_2:=h_2^{-1}$ and change the total area distortion problem into the Dirichlet energy problem of $h_{2}$
\begin{equation}
\label{totaldis:north}
\int_{\hat{\Omega}_0}\frac{|\nabla_{\hat{\mathbf{x}}}h_2(\hat{\mathbf{x}})|_{\delta}^2}{J(\hat{\mathbf{x}}, h_2)}d\hat{x}^2\wedge d\hat{x}^1
=\int_{\mathbb{D}}|\nabla_{\hat{\mathbf{y}}}g_2(\hat{\mathbf{y}})|_{\delta}^2 d\hat{y}^2\wedge d\hat{y}^1.
\end{equation}
From \eqref{northNinv}, \eqref{Totaldis:north} and \eqref{totaldis:north} follows the second term in \eqref{Totaldis1} to be
\begin{align}\label{totaldis:north1}
    \int_{\overline{\mathbb{R}}^2\backslash\Omega_0}\frac{|\nabla_{\mathbf{x}} h_{N}(\mathbf{x})|_\rho^{2}}{J(\mathbf{x},  h_{N})_\rho}dx^1\wedge dx^2
    =\int_{\mathbb{D}}|\nabla_{\hat{\mathbf{y}}}g_2(\hat{\mathbf{y}})|_{\delta}^2 d\hat{y}^2\wedge d\hat{y}^1.
\end{align}
Here, $g_2:=h_2^{-1}$ defines the second term of $ g$ in \eqref{invh}.

Finally, we derive the total area distortion of $ h$ in \eqref{finite2} by summing up the southern and northern parts in \eqref{totaldis:south} and \eqref{totaldis:north1}.
\end{proof}

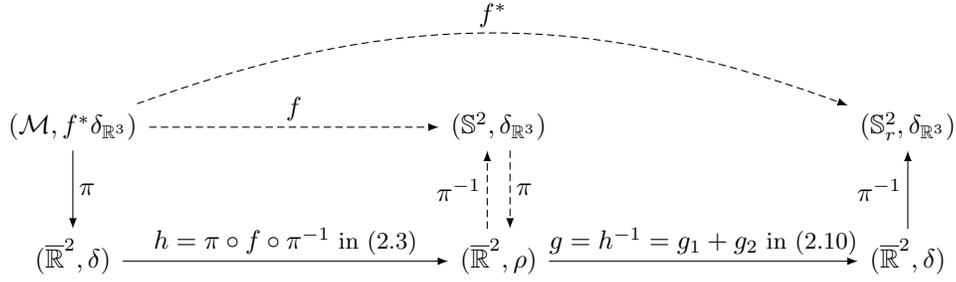
\begin{figure}[h]
\center
\begin{tikzcd}[row sep=3em, column sep=11em]
(\mathcal{M}, f^*\delta_{\mathbb{R}^3}) \arrow{d}{\pi} \arrow[dashed]{r}{f} \arrow[dashed, bend left=20]{rr}{f^*}& (\mathbb{S}^2,  \delta_{\mathbb{R}^3}) \arrow[xshift=1ex, dashed]{d}{\pi}&(\mathbb{S}^2_r, \delta_{\mathbb{R}^3})\\
 (\overline{\mathbb{R}}^2, \delta)\arrow{r}{h=\pi\circ f\circ\pi^{-1} \mbox{\small \ in \eqref{eq:h}}}& (\overline{\mathbb{R}}^2, \rho) \arrow[xshift=-1ex ,dashed]{u}{\pi^{-1}}\arrow{r}{{g = h^{-1} = g_1 + g_2}   \mbox{\small \ in \eqref{invh}}}&(\overline{\mathbb{R}}^2, \delta) \arrow{u}{\pi^{-1}}
\end{tikzcd}
\caption{An illustration for the whole schemetical diagram with $\mathrm{area}(\mathcal{M})=\mathrm{area}(\mathbb{S}^2_r)$, where $\mathbb{S}^2_r$ is a sphere with radius $r$ in $\mathbb{R}^3$.}
\label{Fig1}
\end{figure}

The relationships of $f$, $\pi$, $h$ and $g$ defined in \eqref{eq:h} and \eqref{invh}, respectively, are illustrated in Figure~\ref{Fig1}. 

\begin{theorem}[Equiareal Map]
\label{Equiareal}
Suppose $ h$ is a conformal mapping in $W^{1, 2}(\overline{\mathbb{R}}^2, \overline{\mathbb{R}}^2)$ with $\det\nabla  h(\mathbf{x})>0$ as in Remark~\ref{1} (i). 
There exists the conformal mapping $ g:=g_1+g_2$ in $W^{1, 2}( h(\overline{\mathbb{R}}^2), \overline{\mathbb{R}}^2)$ determined by solving the Euler--Lagrange equation
\begin{align}
    \Delta_{\mathbf{y}}g_1(\mathbf{y})=0,\quad\Delta_{\hat{\mathbf{y}}}g_2(\hat{\mathbf{y}})=0.\label{confh}
\end{align}
Furthermore, the map $ g$ composed with $ h$ is an equiareal (area-preserving) map, i.e.,
\begin{equation}
\begin{cases}
 g\circ h: (\overline{\mathbb{R}}^2, \delta)\to ( g(\overline{\mathbb{R}}^2), \delta),\nonumber\\
\det\nabla g( h(\mathbf{x}))=1,
\end{cases}
\end{equation} 
where $\delta$ is the Euclidean metric in $\mathbb{R}^2$.
\end{theorem}
\begin{proof}
Since $ h$ is a conformal mapping in $W^{1, 2}(\overline{\mathbb{R}}^2, \overline{\mathbb{R}}^2)$ with $\det\nabla  h(\mathbf{x})>0$, we know that $ g$ exists in $W^{1, 2}(\overline{\mathbb{R}}^2, \overline{\mathbb{R}}^2)$ via Theorem~\ref{invconformal}. According to the St\"{o}ilow factorization theorem \cite[chap. 5]{Astala2008}, we have
\begin{align}
     h(\mathbf{x})=
    \begin{cases}
    \phi\circ g_1^{-1}(\mathbf{x}),\quad \mbox{if $\mathbf{x}\in\Omega_0$},\\
    \psi\circ g_2^{-1}(\hat{\mathbf{x}}),\quad \mbox{if $\hat{\mathbf{x}}:=\frac{\mathbf{x}}{|\mathbf{x}|^2}\in\hat{\Omega}_0$},
    \end{cases}
\end{align}
where $\phi$ and $\psi$ are holomorphic functions in $W^{1,2}(h_1(\Omega_0), \overline{\mathbb{R}}^2)$ and $W^{1,2}(h_2(\hat{\Omega}_0), \overline{\mathbb{R}}^2)$, respectively.

In order to determine $ g$ as a conformal mapping in $W^{1, 2}( h(\overline{\mathbb{R}}^2), \overline{\mathbb{R}}^2)$ as in \cite{LHLY}, the equality \eqref{finite2} of Theorem~\ref{invconformal} shows that the minimization of the total area distortion is equivalent to the minimization of the Dirichlet energy of $ g$. Moreover, the Euler--Lagrange equations to the Dirichlet energy of $ g$ are given by \eqref{confh}.

Let $\mathbf{y}= h(\mathbf{x})$. Then by \eqref{recfmap} and \eqref{invh}, we have
\begin{align}
\label{equiareal}
\det\nabla g( h(\mathbf{x}))
&=\det\nabla_\mathbf{y} g(\mathbf{y})\det\nabla_\mathbf{x} h(\mathbf{x})\nonumber\\
&=\frac{1}{\det\nabla_\mathbf{x} h(\mathbf{x})}\det\nabla_\mathbf{x} h(\mathbf{x})=1.
\end{align}
Thus, the Jacobian of $ g\circ h$ is equal to one, i.e., $ g\circ h$ is equiareal.
\end{proof}

For study of the discrete version of the stretch energy functional in Section~\ref{sec:3}, with $J(\mathbf{y},  g)=\frac{1}{J(\mathbf{x},   h)}$ for $ g= h^{-1}$, we have the following corollary.
\begin{corollary}\label{Contdist}
The total area distortion in \eqref{finite2} can be  expressed by
\begin{align}
\int_{\overline{\mathbb{R}}^2}\frac{|\nabla h(\mathbf{x})|_\rho^{2}}{J(\mathbf{x},  h)_\rho}dx^1\wedge dx^2
&=\int_{\mathbb{D}}|\nabla_{\mathbf{y}}g_1(\mathbf{y})|_{\delta}^2dy^1\wedge dy^2+\int_{\mathbb{D}}|\nabla_{\hat{\mathbf{y}}}g_2(\hat{\mathbf{y}})|_{\delta}^2 d\hat{y}^2\wedge d\hat{y}^1 \nonumber\\
&=\int_{\Omega_0}\frac{|\nabla_{\mathbf{y}}g_1(\mathbf{y})|_{\delta}^2}{J(\mathbf{y},g_1)}dx^1\wedge dx^2+\int_{\hat{\Omega}_0}\frac{|\nabla_{\hat{\mathbf{y}}}g_2(\hat{\mathbf{y}})|_{\delta}^2}{J(\hat{\mathbf{y}}, g_2)} d\hat{x}^2\wedge d\hat{x}^1\label{finite1}.
\end{align}
\end{corollary}

\section{Convergence of SEM for spherical equiareal parameterizations} \label{sec:3}

\subsection{Stretch energy functional and SEM algorithm}
A discrete model for a closed genus-zero surface $\mathcal{M}$ can be represented by a closed triangular meshes composed of vertices, edges and faces denoted by $\mathbb{V}(\mathcal{M})$, $\mathbb{E}(\mathcal{M})$ and $\mathbb{F}(\mathcal{M})$, respectively. 
Let $\mathbb{V}(\mathcal{M})=\{\mathrm{v}_t\}_{t=1}^n$. A piecewise affine mapping $f:\mathbb{F}(\mathcal{M})\to\mathbb{R}^3$ on all vertices $\mathrm{v}_t\in\mathbb{V}(\mathcal{M})$ can be expressed as an $n\times3$ matrix
\begin{align*}
\mathbf{f} \equiv [f(\mathrm{v}_1), \cdots,  f(\mathrm{v}_n)]^{\top} \equiv [\mathrm{f}_1, \cdots, \mathrm{f}_n]^{\top}\in\mathbb{R}^{n\times3}.
\end{align*}
Here, $f$ is said to be the induced function by $\f$, and $\f$ is the inducing vector for $f$. 

The original SEM \cite{yueh2019novel} aims to compute a spherical equiareal parameterization of a genus-zero closed surface $\M$ by minimizing the stretch energy functional $E_S(f) =\frac{1}{2} \mathrm{Tr}(\f^\top L(f) \f)$. However, based on the theoretical foundation of \eqref{finite1} and the idea of the uniformization theorem that supports the legitimacy of the definition of the stretch energy, we should modify the stretch energy functional and the associated modified Laplacian matrix as follows.

\begin{subequations} \label{17}
\begin{align}\label{17a-1}
E_S(f)=\frac{1}{2} \left[\mathrm{Tr}(\h_{1}^{\mathrm{H}} L_1(f) \h_{1})+\mathrm{Tr}(\h_{2}^{\mathrm{H}} L_2(f) \h_{2})\right],
\end{align}
where $\h_{1}$ and $\h_{2}$ are the inducing vectors of functions $h_{1}$ and $h_{2}$ as in \eqref{recfmap}, and
$L_s(f)=[w_{ij}^s(f)]_{i, j=1}^{n_s+m_s}$ ($n_s$ and $m_s$, $s=1, 2$ are defined in \eqref{eq:idx_sets_I_Bprime_B} below  for details) with
\begin{align}\label{17b}
    w_{ij}^s(f)=
    \begin{cases}
    -\frac{1}{2} \left(\frac{\cot\theta_{ij}^l(f)}{\sigma_{f^{-1}}([\mathrm{v}_i, \mathrm{v}_l, \mathrm{v}_j])} + \frac{\cot \theta_{ij}^{r}(f)}{\sigma_{f^{-1}}([\mathrm{v}_j, \mathrm{v}_r, \mathrm{v}_i])} \right), & \text{if $[\mathrm{v}_i,\mathrm{v}_j]\in\mathbb{E}(f^{-1}(\pi^{-1}(\mathbb{D})))$, for $s=1$}, \\
    & \quad\text{$[\mathrm{v}_i,\mathrm{v}_j]\in\mathbb{E}(f^{-1}(\hat{\pi}^{-1}(\mathbb{D})))$, for $s=2$},\\
    -\sum_{\ell\neq i} w_{i\ell}^s(f), & \text{if $i=j$}, \\
    0, & \text{otherwise},
    \end{cases}
\end{align}
\end{subequations}
in which $\hat{\pi}(\cdot)=\pi(\cdot)/|\pi(\cdot)|^2$,  $\sigma_{f^{-1}}([\mathrm{v}_i, \mathrm{v}_k, \mathrm{v}_j])=\frac{|[\mathrm{v}_i,\mathrm{v}_k,\mathrm{v}_j]|}{|f([\mathrm{v}_i,\mathrm{v}_k,\mathrm{v}_j])|}$ is the stretch factor (local area ratio) of $f$ on the triangle $[\mathrm{v}_i,\mathrm{v}_k,\mathrm{v}_j]$, for $k= l, r$, $\theta_{ij}^l(f)$ and $\theta_{ij}^r(f)$ are the angles opposite to $[\mathrm{f}_i, \mathrm{f}_j]$ in the adjacent triangles. 
Furthermore, for the boundary edge, i.e., if $[\mathrm{v}_i, \mathrm{v}_j]\in\mathbb{E}(f^{-1}(\pi^{-1}(\partial \mathbb{D})))$ or $\mathbb{E}(f^{-1}(\hat{\pi}^{-1}(\partial \mathbb{D})))$, we set the term involved $\mathrm{v}_r$ $(s=1)$ or $\mathrm{v}_l$ $(s=2)$ in $w_{ij}^s(f)$ to be zero.




\begin{remark}
\begin{enumerate}
    \item[(i)] For any $\Delta \equiv [\mathrm{v}_i, \mathrm{v}_j, \mathrm{v}_k]\in\mathbb{F}(\mathcal{M})$, the stretch factor $\sigma_{f^{-1}}(\Delta)$ in \cite{yueh2019novel} has been proven to be equal to the Jacobian $J_{f^{-1}}|_{f(\Delta)}$. Since we consider $f:\mathcal{M}\to\mathbb{S}^2$ as a conformal mapping, a triangle $f(\Delta)\in \mathbb{F}(f(\mathcal{M}))$ is stereographically projected onto a triangle $h(\Delta_{\pi})\equiv [\mathrm{h}_{i}, \mathrm{h}_{j}, \mathrm{h}_{k}]$ with $\mathrm{h}_{t} \equiv h(\pi(\mathrm{v}_t))\in\mathbb{R}^2$  for $t\in\{i, j, k\}$ and $\Delta_{\pi} \equiv [\pi(\mathrm{v}_i), \pi(\mathrm{v}_j), \pi(\mathrm{v}_k)]$.
    
    Let ``$\mathrm{n}$" be a unit vector normal of the face $[\mathrm{h}_{i}, \mathrm{h}_{j}, \mathrm{h}_{k}]$,  we define the normal vectors on edges of a triangle $[\mathrm{h}_{i}, \mathrm{h}_{j}, \mathrm{h}_{k}]$ as follows, 
    \begin{subequations}\label{balance}
    \begin{align}
           \mathrm{n}_i=\mathrm{n}\times(\mathrm{h}_{k}-\mathrm{h}_{j}),\ 
            \mathrm{n}_j=\mathrm{n}\times(\mathrm{h}_{i}-\mathrm{h}_{k}),\ 
            \mathrm{n}_k=\mathrm{n}\times(\mathrm{h}_{j}-\mathrm{h}_{i}),
    \end{align}
    with
    \begin{align} 
    \mathrm{n}_i+\mathrm{n}_j+\mathrm{n}_k=0.
    \end{align}
    \end{subequations}
    For any point $\mathrm{h}\in \mathrm{Int}\left([\mathrm{h}_{i}, \mathrm{h}_{j}, \mathrm{h}_{k}]\right)$, it can be expressed in the barycentric coordinate,
    \begin{align}
        \mathrm{h}=\lambda_i\mathrm{h}_{ i}+\lambda_j\mathrm{h}_{ j}+\lambda_k\mathrm{h}_{ k}
    \end{align}
     where $\lambda_i=\frac{( \mathrm{h}-\mathrm{h}_{ k})^{\top} \mathrm{n}_i}{2|h(\Delta_{\pi})|}$, $\lambda_j=\frac{( \mathrm{h}-\mathrm{h}_{ i})^{\top} \mathrm{n}_j}{2|h(\Delta_{\pi})|}$, $\lambda_k=\frac{( \mathrm{h}-\mathrm{h}_{ j})^{\top} \mathrm{n}_k}{2|h(\Delta_{\pi})|}$.\\
     
     Similarily, we express $g(\mathrm{h})$ in the barycentric coordinate,
    \begin{align}
    g(\mathrm{h})=\frac{1}{2\left|h(\Delta_{\pi})\right|}\left(\mathrm{g}_i( \mathrm{h}-\mathrm{h}_{ k})^{\top}\mathrm{n}_i+\mathrm{g}_j(\mathrm{h}-\mathrm{h}_{ i})^{\top}\mathrm{n}_j+\mathrm{g}_k(\mathrm{h}-\mathrm{h}_{ j})^{\top} \mathrm{n}_k\right).
    \end{align}
    Taking the gradient of $\mathrm{h}$, we have
    \begin{align}\label{disgrad}
        \nabla_{\mathrm{h}}g=\frac{1}{2\left|h(\Delta_{\pi})\right|}\sum_{t\in\{i, j, k\}}\mathrm{g}_t\mathrm{n}_t.
    \end{align}
    The discrete version of Dirichlet energy $g$ of \eqref{finite2} restricted on $h(\Delta_{\pi})$ is given by
    \begin{align}\label{disDiri}
        E(g)\rvert_{h(\Delta_{\pi})}=\frac{1}{2}\langle \nabla_{\mathrm{h}}g, \nabla_{\mathrm{h}}g\rangle\left|h(\Delta_{\pi})\right|.
    \end{align}
    
    By \eqref{recfmap} and \eqref{invh}, we take $g=h^{-1}$ be conformal mapping and check $\langle \nabla_{\mathrm{h}}g,\nabla_{\mathrm{h}}g\rangle\left|h(\Delta_{\pi})\right|$ as follows,
    \begin{align}
       & \langle \nabla h^{-1}, \nabla h^{-1}\rangle\left|h(\Delta_{\pi})\right|
        =-\frac{1}{2}\sum_{\substack{s\neq t\\ s, t\in\{i, j, k\}}}\frac{\langle \mathrm{n}_s, \mathrm{n}_t\rangle}{2\left|h(\Delta_{\pi})\right|}(\mathrm{g}_s-\mathrm{g}_t)^2\hspace{6mm} (\text{by  \eqref{balance} and \eqref{disgrad}})\nonumber\\
        &\quad=-\frac{1}{2}\sum_{\substack{s\neq t\\ s,t\in\{i, j, k\}}}-\cot\theta_{st}^{k}(h)[J(\mathrm{g}, h)(\mathrm{h}_{ s}-\mathrm{h}_{ t})^2]\nonumber\\
        &\quad\quad\quad\quad\quad\quad\quad\quad\quad\quad \left(\text{because of   $(\mathrm{g}_s-\mathrm{g}_t)=\sqrt{J(\mathrm{g}, h)}(\mathrm{h}_{ s}-\mathrm{h}_{ t})$}\right)\nonumber\\
        &\quad=-\frac{1}{2}\sum_{\substack{s\neq t\\ s,t\in\{i, j, k\}}}\frac{-\cot\theta_{st}^{k}(h)}{J(\mathrm{h}, g)}(\mathrm{h}_{ s}-\mathrm{h}_{ t})^2 \hspace{10mm} \left(J(\mathrm{h}, g)=\frac{1}{J(\mathrm{g}, h)}\right)\nonumber\\
        &\quad=-\frac{1}{2}\sum_{\substack{s\neq t\\ s,t\in\{i, j, k\}}}\frac{-\cot\theta_{st}^{k}(h)}{\left|\Delta_{\pi}\right|/\left|h(\Delta_{\pi})\right|}(\mathrm{h}_{ s}-\mathrm{h}_{ t})^2 \hspace{6mm} \left(J_{h^{-1}}\rvert_{h(\Delta_{\pi})}=\frac{\left|\Delta_{\pi}\right|}{\left|h(\Delta_{\pi})\right|}, \text{by \cite{yueh2019novel} }\right)\nonumber\\
        &\quad=-\frac{1}{2}\sum_{\substack{s\neq t\\ s,t\in\{i, j, k\}}}\frac{-\cot\theta_{st}^{k}(f)}{\left|\Delta\right|/\left|f(\Delta)\right|}(\mathrm{h}_{ s}-\mathrm{h}_{ t})^2.\label{lasteq}
    \end{align}
    The last equality \eqref{lasteq} holds, because $f$ and $h$ are different up to a  stereographic projection. Therefore, we derive the equations \eqref{disDiri} and \eqref{lasteq} that show the modified SEM \eqref{17a-1} according to the mathematical foundation \eqref{finite1} of Corollary~\ref{Contdist}.

    \item[(ii)] In light of \eqref{invh} and Theorem~\ref{Equiareal}, if $ g(y)=g_1(\mathbf{y})+g_2(\hat{\mathbf{y}})$ solves the minimization of the Dirichlet energy in \eqref{finite2}, then $ g\circ h$ is an equiareal map as we expected.
\end{enumerate}
\end{remark}

Since the energy functional in \eqref{17a-1} is highly nonlinear in $\h=\h_{1}+\h_{2}$, solving the SEM problem \eqref{17} is much more complicate than that of the CEM in \cite{LHLY}. Then We propose the following strategy to solve \eqref{17}, which is a slightly modified version of the SEM in \cite{YuLL19,yueh2019novel}.

\begin{algorithm}[]
\caption{Modified SEM for spherical equiareal parameterizations}
\label{alg:SEM_fix_bnd_k}
\begin{algorithmic}[1]
\REQUIRE A closed genus-zero surface $\M$, a tolerance $\varepsilon$, a radius $r$ (e.g., $\varepsilon=10^{-6}$, $r=1.1$). 
\ENSURE A spherical equiareal parameterization $\f^{(k)}$.
\STATE Set $n=\#\V(\M)$ and $\delta = \infty$.
\STATE Compute a spherical conformal parameterization $\f$ using CEM algorithm in \cite{YuLL19} for satisfying the condition of Theorem~\ref{Equiareal}.
\STATE Compute $\h_\ell=\frac{\f_{\ell,1}}{1-\f_{\ell,3}}+\i\frac{\f_{\ell,2}}{1-\f_{\ell,3}}$, $\ell=1, \ldots, n$. 
\STATE \% Following for-loop is used to fix the indices of $\mathtt{B}_s$ and $\mathtt{I}_s$.
\FOR{$k=1,2,3$ } 
\FOR{$s=1,2$}
\STATE Update $L \gets L_s(f)$, where $L_s(f)$ is defined as in \eqref{17b}.
\STATE Update $\h\gets\diag(|\h|^{-2})\h$.  
\STATE Let $\I_s=\{i\mid|\h_i|<r\}$, $\B_s^{\prime}=\{1, \ldots,n\}\backslash\I_s$, and $\mathtt{B}_s = \{ j \mid [v_i, v_j ] \in \mathbb{E}(\mathcal{S}), i \in \mathtt{I}_s, j \in \mathtt{B}_s^{\prime}\}$.  
\STATE Update $\h$ by solving
$[L_s]_{\I_s,\I_s}\h_{\I_s} = -[L_s]_{\I_s,\B_s} \h_{\B_s}$.  
\STATE Compute $\f_\ell=\Pi_{\bS^2}^{-1}(\h_\ell)$, $\ell=1, \ldots, n$.  
\ENDFOR
\ENDFOR 
\STATE Set $k = 0$, $\f_{1}^{(0)} = \f_{2}^{(0)} = \f$ and $\h_{1}^{(0)} = \h_{2}^{(0)} = \h$.
\STATE Compute $L_s(f)$ in \eqref{17b}. 
\STATE Compute $\mathcal{E}^{(0)} = E_S(f)$.
\WHILE{$\delta>\varepsilon$} 
    \STATE Set $ k = k + 1$.  
    \STATE Compute $\h_{1}^{(k)} = \diag(|\h_{2}^{(k-1)}|^{-2}) \h_{2}^{(k-1)}$. \ \% Southern hemispherical iteration \label{alg:south_inversion_north}
    \STATE Update $\h_{1}^{(k)}$ by solving $[L_{1}(f^{(k-1)})]_{\I_1,\I_1} \h_{\I_1}^{(k)} = -[L_{1}(f^{(k-1)})]_{\I_1,\B_1}  \h_{\B_1}^{(k)}$. \label{alg:south_solving_LS}
    \STATE Compute $\f_\ell^{(k)}=\Pi_{\bS^2}^{-1}({\h}_\ell^{(k)})$, $\ell=1, \ldots, n$, and $L_{1}(f^{(k)})$ in \eqref{17b}.
    \label{algSEM:21}
    \STATE Compute $\h_{2}^{(k)} = \diag(|\h_{1}^{(k)}|^{-2}) \h_{1}^{(k)}$.  \quad \% Northern hemispherical iteration \label{alg:north_inversion_south}
    \STATE Update $\h_{2}^{(k)}$ by solving $[L_{2}(f^{(k-1)})]_{\I_2,\I_2} \h_{\I_2}^{(k)} = -[L_{2}(f^{(k-1)})]_{\I_2,\B_2} \h_{\B_2}^{(k)}$. \label{alg:north_solving_LS}
\STATE Compute $\f_\ell^{(k)}=\Pi_{\bS^2}^{-1}(\h_\ell^{(k)})$, $\ell=1, \ldots, n$, and $L_2(f^{(k)})$ in \eqref{17b}. \label{algSEM:inv_ste}
\STATE Compute $\mathcal{E}^{(k)} = E_S(f^{(k)})$ and $\delta= \mathcal{E}^{(k)}-\mathcal{E}^{(k-1)}$.   
\ENDWHILE
\RETURN the spherical equiareal parameterization $\f^{(k)}$.
\end{algorithmic}
\end{algorithm}

In view of \eqref{finite1} and \eqref{lasteq}, the energy functional \eqref{17a-1} can be minimized by alternatively solving $\mathbf{h}_{1}$ and $\mathbf{h}_{2}$ on the complex plane in association with modified Laplacian matrices in \eqref{17b} corresponding to the southern and the inversion of the northern hemispheres, respectively, as in steps \ref{alg:south_inversion_north}--\ref{alg:north_solving_LS} of Algorithm \ref{alg:SEM_fix_bnd_k}.  

Given a radius $\rho\gtrsim 1$, we define
\begin{subequations} \label{eq:idx_sets_I_Bprime_B}
\begin{align}
\mathtt{I}_s &=\{i\mid\lvert {h}_i^{(s)}\rvert < \rho\}, \ \mathtt{B}_s^{\prime}=\{1, \dots,n\}\backslash\mathtt{I}_s, \\ 
\mathtt{B}_s &= \{ j \mid [\mathrm{v}_i, \mathrm{v}_j ] \in \mathbb{E}(\mathcal{M}), i \in \mathtt{I}_s, j \in \mathtt{B}_s^{\prime}\}, 
\end{align}
\end{subequations}
where $n=\#\V(\M)$,  $n_s = \# \mathtt{I}_{s}$ and $m_s = \# \mathtt{B}_s$, and the indices $s=1, 2,$ indicate the southern and northern hemispheres, respectively. Let $f^{(k)}$ be the induced function of $\f^{(k)}$ at the $k$th step of Algorithm~\ref{alg:SEM_fix_bnd_k}, and let
\begin{subequations} \label{22}
\begin{align} 
L_s^{(k)} &= [L_s(f^{(k)})]_{\I_s,\I_s}, \quad
B_s^{(k)} = [L_s(f^{(k)})]_{\I_s,\B_s}, \label{22a} \\
P_1 &= [I_n]_{\B_{2},\I_1}, ~ P_2 = [I_n]_{\B_{1},\I_2}, 
\label{22b}
\end{align} 
and
\begin{equation} \label{22c}
A_s^{(k)} = -P_s \left(L_s^{(k)}\right)^{-1} B_s^{(k)} \equiv P_s \widehat{A}_{s}^{(k)}, ~ s=1,2.
\end{equation}
\end{subequations}
The $k$th southern and northern hemispherical iterations of steps \ref{alg:south_inversion_north}-\ref{algSEM:21} and \ref{alg:north_inversion_south}-\ref{algSEM:inv_ste} can be written as
\begin{subequations} \label{24}
\begin{alignat}{3}
\h_{\B_1}^{(k)} &= \diag(|P_2 \h_{\I_2}^{(k-1)}|)^{-2} P_2 \h_{\I_2}^{(k-1)}, &\quad
\h_{\I_{1}}^{(k)} &= \widehat{A}_{1}^{(k-1)} \h_{\B_1}^{(k)},  
\label{24a} \\
\h_{\B_2}^{(k)} &= \diag(|P_1 \h_{\I_1}^{(k)}|)^{-2} P_1 \h_{\I_1}^{(k)},  & \quad 
\h_{\I_{2}}^{(k)} &= \widehat{A}_{2}^{(k-1)} \h_{\B_2}^{(k)}.  
\label{24b}
\end{alignat}
\end{subequations}

\subsection{R-linear convergence of Algorithm~\ref{alg:SEM_fix_bnd_k}}
To analyze the convergence of the iterations in \eqref{24}, let 
\begin{align}  
\beps_{\B_s}^{(k+1)} &= \h_{\B_s}^{(k+1)} - \h_{\B_s}^{(k)}   
\label{32}  
\end{align}
for $s = 1, 2$.
From \eqref{24} follows that,
\begin{subequations} \label{eq:eps_1_2}
\begin{align}
\beps_{\B_1}^{(k+1)} &= \h_{\B_1}^{(k+1)}-\h_{\B_1}^{(k)} \nonumber \\
  &= \diag(|A_2^{(k-1)}\h_{\B_2}^{(k)}|)^{-2} A_2^{(k-1)} \h_{\B_2}^{(k)} - \diag(|A_2^{(k-2)}\h_{\B_2}^{(k-1)}|)^{-2} A_2^{(k-2)} \h_{\B_2}^{(k-1)} \nonumber \\
 &= \mbox{conj}(\Gamma_2^{(k)} \left( A_2^{(k-2)} \h_{\B_2}^{(k-1)} - A_2^{(k-1)} \h_{\B_2}^{(k)} \right) ) \nonumber \\
&= \mbox{conj}(\Gamma_2^{(k)} \left( -A_2^{(k-1)} \beps_2^{(k)}  + (A_2^{(k-2)} - A_2^{(k-1)}) \h_{\B_2}^{(k-1)} \right)) \label{34-1}
\end{align} 
and
\begin{align}
 \beps_{\B_2}^{(k+1)} &= \h_{\B_2}^{(k+1)}-\h_{\B_2}^{(k)} \nonumber \\
 &= \diag(|A_1^{(k)}\h_{\B_1}^{(k+1)}|)^{-2} A_1^{(k)} \h_{\B_1}^{(k+1)} - \diag(|A_1^{(k-1)}\h_{\B_1}^{(k)}|)^{-2} A_1^{(k-1)} \h_{\B_1}^{(k)} \nonumber \\ 
&= \mbox{conj}(\Gamma_1^{(k+1)} \left( -A_1^{(k)} \beps_{\B_1}^{(k+1)}  + (A_1^{(k-1)} - A_1^{(k)}) \h_{\B_1}^{(k)} \right)), \label{34}
\end{align} 
\end{subequations}
where
\begin{align*}
      \Gamma_s^{(k)} & = \diag\left( \frac{1}{\left( \mathbf{e}_i^{\top} A_s^{(k-2)} \h_{\B_s}^{(k-1)}\right)\mathbf{e}_i^{\top} A_s^{(k-1)} \h_{\B_s}^{(k)}}\right)_{i=1}^{n_{\B_{s+1}}}.
\end{align*}

The essential term in \eqref{eq:eps_1_2} is $(A_s^{(k)} - A_s^{(k-1)}) \h_{\B_s}^{(k)} = P_s(\widehat{A}_s^{(k)} - \widehat{A}_s^{(k-1)}) \h_{\B_s}^{(k)}$. From \eqref{22c}, this term can be reformulated as 
\begin{subequations}
\begin{align}
& (\widehat{A}_s^{(k)} - \widehat{A}_s^{(k-1)}) \h_{\B_s}^{(k)} = \left( \left( L_s^{(k-1)} \right)^{-1} B_s^{(k-1)} -\left( L_s^{(k)} \right)^{-1} B_s^{(k)} \right) \h_{\B_s}^{(k)} \nonumber \\
=& \left( \left( L_s^{(k-1)} \right)^{-1} -\left( L_s^{(k)} \right)^{-1}  \right) B_s^{(k-1)} \h_{\B_s}^{(k)} - \left( L_s^{(k)} \right)^{-1} \left( B_s^{(k)} - B_s^{(k-1)} \right) \h_{\B_s}^{(k)} \nonumber \\
=& \left( L_s^{(k)} \right)^{-1} \left( L_s^{(k)} -  L_s^{(k-1)}  \right) \left( L_s^{(k-1)} \right)^{-1}  B_s^{(k-1)} \h_{\B_s}^{(k)} - \left( L_s^{(k)} \right)^{-1} \left( B_s^{(k)} - B_s^{(k-1)} \right) \h_{\B_s}^{(k)} \nonumber \\
=& \left( L_s^{(k)} \right)^{-1} \begin{bmatrix} 
      L_s^{(k)} -  L_s^{(k-1)} & B_s^{(k)} - B_s^{(k-1)}
\end{bmatrix} \mathbf{g}_s^{(k)} \nonumber \\
=& \left( L_s^{(k)} \right)^{-1} \left( \begin{bmatrix} 
      L_s^{(k)} & B_s^{(k)} 
\end{bmatrix} - \begin{bmatrix} 
      L_s^{(k-1)} & B_s^{(k-1)} 
\end{bmatrix} \right) \mathbf{g}_s^{(k)}, \label{36}
\end{align}
where
\begin{align}
    \mathbf{g}_s^{(k)} = \begin{bmatrix}
           \left( L_s^{(k-1)} \right)^{-1}  B_s^{(k-1)} \h_{\B_s}^{(k)} \\ - \h_{\B_s}^{(k)} 
      \end{bmatrix} \equiv [g_{\ell}]. \label{eq:h_hat_k}  
\end{align}
\end{subequations}
Therefore, we can focus on the analysis of $([L_s^{(k)}, B_s^{(k)}] - [L_s^{(k-1)}, B_s^{(k-1)}])\mathbf{g}_s^{(k)}$.
From \eqref{17b}, we have
\begin{align}
     \left( \begin{bmatrix} 
      L_s^{(k)} & B_s^{(k)} 
\end{bmatrix} \right. & \left. - \begin{bmatrix} 
      L_s^{(k-1)} & B_s^{(k-1)} 
\end{bmatrix} \right) \mathbf{g}_s^{(k)} = \left[\widehat{w}_{ii}^{(k)} g_i + \sum_{j \in \N(i) } \widehat{w}_{ij}^{(k)} g_j \right]_{i=1}^{n_1} \nonumber \\
     &= \left[ \sum_{j \in \N(i) } \widehat{w}_{ij}^{(k)} g_j - \sum_{j \in \N(i) } \widehat{w}_{ij}^{(k)} g_i \right]_{i=1}^{n_1} \nonumber \\
     &= \left[ \sum_{j \in \N(i) } \widehat{w}_{ij}^{(k)} (g_j - g_i) \right]_{i=1}^{n_1} \label{eq:Lk_m_Lkm1_h}
\end{align}
with $\widehat{w}_{ij}^{(k)} = w_{ij}^s(f^{(k)}) - w_{ij}^s(f^{(k-1)})$. Now, we give an analysis about $\widehat{w}_{ij}^{(k)}$ in Lemma~\ref{lem:w_ij}.

\begin{lemma} \label{lem:w_ij}
      For $[\mathrm{v}_i,\mathrm{v}_j]\in\E(\M)$, $\Delta_l \equiv [\mathrm{v}_i,\mathrm{v}_l,\mathrm{v}_j]$ and $\Delta_r\equiv[\mathrm{v}_j,\mathrm{v}_r,\mathrm{v}_i] \in \F(\M)$,
      let $w_{ij}^s(f^{(k)})$ be the $(i,j)$-th element of the Laplacian matrix $L_{s}(f^{(k)})$  in \eqref{17b} and
      \begin{align*} 
      \varepsilon_t^{(k)} = \mathrm{h}_t^{(k)} - \mathrm{h}_t^{(k-1)}, \quad t = 1, \ldots, n
      \end{align*} 
      with 
      \begin{align}
            \varepsilon_t^{(k)} = | \varepsilon_t^{(k)} | e^{\sqrt{-1}\psi_{t}^{(k)}}. \label{eq:polor_cor_eps}
      \end{align}
      Define
      \begin{align}
       d_{t}^{(k)} = \frac{2 e^{-\sqrt{-1}\psi_t^{(k)}}}{\sqrt{1 + | \mathrm{h}_{t}^{(k)} |^2 } \sqrt{1 + | \mathrm{h}_{t}^{(k-1)} |^2}}. \label{eq:para_dt}
      \end{align}
      Then
      \begin{align}
     w_{ij}(f^{(k)}) - w_{ij}(f^{(k-1)}) 
     = c_i^{(k)} \varepsilon_i^{(k)}  - c_l^{(k)} \varepsilon_l^{(k)}   + c_j^{(k)} \varepsilon_j^{(k)}  - c_r^{(k)} \varepsilon_r^{(k)}, \label{eq:b_hat_ij_2}
\end{align}
where 
\begin{align*}
     c_l^{(k)} 
     &= \frac{ d_l^{(k)}  \| \rf_l^{(k)}-\rf_j^{(k)} + \rf_l^{(k-1)} - \rf_i^{(k-1)} \| \cos \varphi_{l,l}^{(k)}}{4   | \Delta_l |}, \\
     c_i^{(k)} 
     &= \frac{ d_i^{(k)} \| \rf_l^{(k)}-\rf_j^{(k)}  \| \cos \varphi_{l,i}^{(k)}}{4 | \Delta_l | } + \frac{d_i^{(k)} \| \rf_r^{(k-1)} - \rf_j^{(k-1)} \| \cos \varphi_{r,i}^{(k)} }{4 | \Delta_r | }, \\
     c_j^{(k)} 
     &= \frac{d_j^{(k)} \| \rf_l^{(k-1)} - \rf_i^{(k-1)} \| \cos \varphi_{l,j}^{(k)}}{4 | \Delta_l |  } + \frac{d_j^{(k)} \| \rf_r^{(k)}-\rf_i^{(k)}  \| \cos \varphi_{r,j}^{(k)}}{4 | \Delta_r | }, \\
     c_r^{(k)} 
     &= \frac{d_r^{(k)} \| \rf_r^{(k)}-\rf_i^{(k)} + \rf_r^{(k-1)} - \rf_j^{(k-1)} \| \cos \varphi_{r,r}^{(k)}}{4 | \Delta_r | }
\end{align*}
with $\varphi_{l,l}$, $\varphi_{l,i}$, $\varphi_{l,j}$ being the angles between $(\rf_l^{(k)} - \rf_{l}^{(k-1)})$ and $(\rf_l^{(k)}-\rf_j^{(k)} + \rf_l^{(k-1)} - \rf_i^{(k-1)})$, $(\rf_i^{(k)} - \rf_i^{(k-1)})$ and $(\rf_l^{(k)}-\rf_j^{(k)})$,  $(\rf_j^{(k)} - \rf_j^{(k-1)})$ and $(\rf_l^{(k-1)} - \rf_i^{(k-1)})$, respectively, and $\varphi_{r,i}$, $\varphi_{r,j}$, $\varphi_{r,r}$ being the angles between $(\rf_i^{(k)} - \rf_i^{(k-1)})$ and $(\rf_r^{(k-1)} - \rf_j^{(k-1)})$, $(\rf_j^{(k)} - \rf_j^{(k-1)})$ and $(\rf_r^{(k)}-\rf_i^{(k)})$, $(\rf_r^{(k)} - \rf_r^{(k-1)})$ and $(\rf_r^{(k)}-\rf_i^{(k)} + \rf_r^{(k-1)} - \rf_j^{(k-1)})$, respectively.
\end{lemma}
\begin{proof}
 For $[\mathrm{v}_i,\mathrm{v}_j]\in\E(\M)$, $\Delta_l \equiv [\mathrm{v}_i,\mathrm{v}_l,\mathrm{v}_j]$ and $\Delta_r\equiv[\mathrm{v}_j,\mathrm{v}_r,\mathrm{v}_i] \in \F(\M)$, from the fact that 
\begin{subequations} \label{37}
\begin{equation} \label{37a}
(\rf_l - \rf_i)^\top(\rf_l-\rf_j) = \|\rf_l - \rf_i\|\|\rf_l-\rf_j\| \cos\theta_{ij}^{l}(f),
\end{equation}
and
\begin{equation} \label{37b}
 2 |f_1(\Delta_l)| = \| (\rf_l - \rf_i)\times (\rf_l-\rf_j) \| = \|\rf_l - \rf_i\|\|\rf_l-\rf_j\| \sin\theta_{ij}^{l}(f),
\end{equation}
\end{subequations}
$w_{ij}^s(f^{(k)})$ in \eqref{17b} with stretch factor $\sigma_{f^{-1}}(\Delta) = |\Delta|/|f(\Delta)|$ can be represented as 
\begin{align}
-w_{ij}^s(f^{(k)})  &= \frac12 \left( \frac{\cos\theta_{ij}^{l}(f^{(k)})}{\sin\theta_{ij}^{l}(f^{(k)})} \frac{|f^{(k)}(\Delta_l)|}{|\Delta_l|} + \frac{\cos\theta_{ij}^{r}(f^{(k)})}{\sin\theta_{ij}^{r}(f^{(k)})}  \frac{|f^{(k)}(\Delta_r)|}{|\Delta_r|} \right) \nonumber \\
%
&= \frac{(\rf_l^{(k)} - \rf_i^{(k)})^\top(\rf_l^{(k)}-\rf_j^{(k)})}{4 |\Delta_l|} + \frac{(\rf_r^{(k)} - \rf_i^{(k)})^\top(\rf_r^{(k)}-\rf_j^{(k)})}{4 |\Delta_r|}. \label{38}
\end{align}
This implies 
\begin{align}
     & w_{ij}^s(f^{(k)}) - w_{ij}^s(f^{(k-1)}) \nonumber \\
     =&\  -\frac{(\rf_l^{(k)} - \rf_i^{(k)})^\top(\rf_l^{(k)}-\rf_j^{(k)}) - (\rf_l^{(k-1)} - \rf_i^{(k-1)})^\top(\rf_l^{(k-1)}-\rf_j^{(k-1)})}{4 |\Delta_l|} \nonumber \\
     -&\ \frac{(\rf_r^{(k)} - \rf_i^{(k)})^\top(\rf_r^{(k)}-\rf_j^{(k)}) - (\rf_r^{(k-1)} - \rf_i^{(k-1)})^\top(\rf_r^{(k-1)}-\rf_j^{(k-1)})}{4 |\Delta_r|}.  \label{eq:err_wij}
\end{align}
Rewrite the molecular of the first and the second term as
\begin{subequations} \label{eq:pf_err_wij_0}
\begin{align}
      &\ (\rf_l^{(k)} - \rf_i^{(k)})^\top(\rf_l^{(k)}-\rf_j^{(k)}) - (\rf_l^{(k-1)} - \rf_i^{(k-1)})^\top(\rf_l^{(k-1)}-\rf_j^{(k-1)}) \nonumber \\
      =&\ (\rf_l^{(k)} - \rf_{l}^{(k-1)} - \rf_i^{(k)} + \rf_i^{(k-1)})^\top(\rf_l^{(k)}-\rf_j^{(k)}) \nonumber \\
      & - (\rf_l^{(k-1)} - \rf_i^{(k-1)})^\top(\rf_l^{(k-1)}- \rf_l^{(k)} + \rf_j^{(k)} - \rf_j^{(k-1)}) \nonumber \\
      =&\ (\rf_l^{(k)} - \rf_{l}^{(k-1)})^{\top} (\rf_l^{(k)}-\rf_j^{(k)} + \rf_l^{(k-1)} - \rf_i^{(k-1)}) \nonumber \\
      & - (\rf_i^{(k)} - \rf_i^{(k-1)})^{\top} (\rf_l^{(k)}-\rf_j^{(k)}) - (\rf_j^{(k)} - \rf_j^{(k-1)})^{\top}(\rf_l^{(k-1)} - \rf_i^{(k-1)}) 
\end{align}
and
\begin{align}
      &\ (\rf_r^{(k)} - \rf_j^{(k)})^\top(\rf_r^{(k)}-\rf_i^{(k)}) - (\rf_r^{(k-1)} - \rf_j^{(k-1)})^\top(\rf_r^{(k-1)}-\rf_i^{(k-1)}) \nonumber \\
      =&\ (\rf_r^{(k)} - \rf_r^{(k-1)})^{\top} (\rf_r^{(k)}-\rf_i^{(k)} + \rf_r^{(k-1)} - \rf_j^{(k-1)}) \nonumber \\
      & - (\rf_j^{(k)} - \rf_j^{(k-1)})^{\top} (\rf_r^{(k)}-\rf_i^{(k)}) - (\rf_i^{(k)} - \rf_i^{(k-1)})^{\top}(\rf_r^{(k-1)} - \rf_j^{(k-1)}). 
\end{align}
\end{subequations}
On the other hand, by the inverse stereographic projection and \eqref{eq:polor_cor_eps}-\eqref{eq:para_dt}, it holds that
\begin{align*}
       \|  \rf_{i}^{(k)}  -\rf_{i}^{(k-1)}  \| &= 2    \left( \sqrt{1 + | \mathrm{h}_{i}^{(k)} |^2 } \sqrt{1 + | \mathrm{h}_{i}^{(k-1)} |^2} \right)^{-1} |  \varepsilon_{i}^{(k)} | = d_i^{(k)} \varepsilon_{i}^{(k)}. 
\end{align*}
This implies that
\begin{align*}
      (\rf_i^{(k)} - \rf_i^{(k-1)})^{\top} (\rf_l^{(k)}-\rf_j^{(k)}) &= \| \rf_i^{(k)} - \rf_i^{(k-1)} \| \left( \| \rf_l^{(k)}-\rf_j^{(k)}  \| \cos \theta_{l,2}^{(k)}\right) \nonumber \\
      &= \left( d_i^{(k)}   \| \rf_l^{(k)}-\rf_j^{(k)}  \| \cos \varphi_{l,i}^{(k)}  \right) \varepsilon_i^{(k)}.
\end{align*}
Therefore, \eqref{eq:pf_err_wij_1} can be represented as
\begin{subequations} \label{eq:pf_err_wij_1}
\begin{align}
     &\ (\rf_l^{(k)} - \rf_i^{(k)})^\top(\rf_l^{(k)}-\rf_j^{(k)}) - (\rf_l^{(k-1)} - \rf_i^{(k-1)})^\top(\rf_l^{(k-1)}-\rf_j^{(k-1)}) \nonumber \\
      =&\ \left( d_l^{(k)}  \| \rf_l^{(k)}-\rf_j^{(k)} + \rf_l^{(k-1)} - \rf_i^{(k-1)} \| \cos \varphi_{l,l}^{(k)}  \right) \varepsilon_{l}^{(k)} \nonumber \\
      & -\left( d_i^{(k)} \| \rf_l^{(k)}-\rf_j^{(k)}  \| \cos \varphi_{l,i}^{(k)} \right) \varepsilon_{i}^{(k)} 
       - \left( d_j^{(k)} \| \rf_l^{(k-1)} - \rf_i^{(k-1)} \| \cos \varphi_{l,j}^{(k)}  \right) \varepsilon_j^{(k)} 
\end{align}
and
\begin{align}
      &\ (\rf_r^{(k)} - \rf_j^{(k)})^\top(\rf_r^{(k)}-\rf_i^{(k)}) - (\rf_r^{(k-1)} - \rf_j^{(k-1)})^\top(\rf_r^{(k-1)}-\rf_i^{(k-1)}) \nonumber \\
      =&\ \left( d_r^{(k)} \| \rf_r^{(k)}-\rf_i^{(k)} + \rf_r^{(k-1)} - \rf_j^{(k-1)} \| \cos \varphi_{r,r}^{(k)}   \right) \varepsilon_{r}^{(k)} \nonumber \\
      & - \left( d_j^{(k)} \| \rf_r^{(k)}-\rf_i^{(k)}  \| \cos \varphi_{r,j}^{(k)}  \right) \varepsilon_{j}^{(k)} 
       - \left( d_i^{(k)} \| \rf_r^{(k-1)} - \rf_j^{(k-1)} \| \cos \varphi_{r,i}^{(k)} \right) \varepsilon_{i}^{(k)}.      
\end{align}
\end{subequations}
Substituting \eqref{eq:pf_err_wij_1} into \eqref{eq:err_wij}, the result in \eqref{eq:b_hat_ij_2} can be obtained. 
\end{proof}

Using the result in Lemma~\ref{lem:w_ij}, we obtain
\begin{align}
     \sum_{j \in \N(i) } \widehat{w}_{ij}^{(k)} (g_j - g_i) &= \sum_{j \in \N(i) } (g_j - g_i) \left( c_i^{(k)} \varepsilon_i^{(k)} - c_l^{(k)} \varepsilon_l^{(k)}   + c_j^{(k)} \varepsilon_j^{(k)}  - c_r^{(k)} \varepsilon_r^{(k)} \right)\nonumber \\
     &= c_i^{(k)} \varepsilon_i^{(k)}  \sum_{j \in \N(i) } (g_j - g_i) + \sum_{j \in \N(i) }  c_j^{(k)} (g_j - g_i) \varepsilon_j^{(k)}  \nonumber \\
     &- \sum_{j \in \N(i), \{r, l\} = \N(i) \cap \N(j) }  (g_j - g_i) \left( c_l^{(k)} \varepsilon_l^{(k)}  + c_r^{(k)} \varepsilon_r^{(k)}  \right) \nonumber \\
      &= \left( c_i^{(k)} \varepsilon_i^{(k)}  \sum_{j \in \N(i) } (g_j - g_i) + \sum_{j \in \N(i)\cap \I_s }  c_j^{(k)} (g_j - g_i) \varepsilon_j^{(k)} \right) \nonumber \\
      &  + \sum_{j \in \N(i)\cap \B_s }  c_j^{(k)} (g_j - g_i) \varepsilon_j^{(k)}  \nonumber \\
     &- \sum_{j \in \N(i), \{r, l\} = \N(i) \cap \N(j) }  (g_j - g_i) \left( c_l^{(k)} \varepsilon_l^{(k)}  + c_r^{(k)} \varepsilon_r^{(k)}  \right) \nonumber \\
     &\equiv {H}_{\I_s}^{(k)}(i,:)  \beps_{\I_s}^{(k)} + {H}_{\B_s}^{(k)}(i,:)  \beps_{\B_s}^{(k)}. \label{eq:mtx_Hk}
\end{align}
Plugging \eqref{eq:mtx_Hk} into \eqref{eq:Lk_m_Lkm1_h}, we have
\begin{align}
     \left( \begin{bmatrix} 
      L_s^{(k)} & B_s^{(k)} 
\end{bmatrix} - \begin{bmatrix} 
      L_s^{(k-1)} & B_s^{(k-1)} 
\end{bmatrix} \right) \mathbf{g}_s^{(k)} &= {H}_{\I_s}^{(k)}  \beps_{\I_s}^{(k)} + {H}_{\B_s}^{(k)}  \beps_{\B_s}^{(k)}\equiv {H}_{s}^{(k)}\beps_{s}^{(k)}. \label{eq:L_k-L_km1_2_Hs}
\end{align}
Combining \eqref{eq:L_k-L_km1_2_Hs} with \eqref{36}, we get following important theorem.

\begin{theorem} \label{thm:tran_Ak_km1_h}
Let
\begin{align}
\beps_s^{(k)}  &\equiv \begin{bmatrix}
\beps_{\I_{s}}^{(k)}\\
\beps_{\B_s}^{(k)}
\end{bmatrix} \equiv \begin{bmatrix}
\h_{\I_{s}}^{(k)}\\
\h_{\B_s}^{(k)}
\end{bmatrix}
-\begin{bmatrix}
\h_{\I_{s}}^{(k-1)}\\
\h_{\B_s}^{(k-1)}
\end{bmatrix}, \ \mathbf{g}_s^{(k)} = \begin{bmatrix}
           \left( L_s^{(k-1)} \right)^{-1}  B_s^{(k-1)} \h_{\B_s}^{(k)} \\ - \h_{\B_s}^{(k)} 
      \end{bmatrix} \label{eq:eta_k}  
\end{align}
for $s = 1, 2$, and $H_{s}^{(k)} \equiv [H_{\I_s}^{(k)}, H_{\B_s}^{(k)}]$ be defined in \eqref{eq:L_k-L_km1_2_Hs}. 
Then
\begin{align}
     (\widehat{A}_s^{(k)} - \widehat{A}_s^{(k-1)}) \h_{\B_s}^{(k)} &= \left( L_s^{(k)} \right)^{-1} H_{s}^{(k)} \beps_s^{(k)}. 
     \label{eq:Ak_Akm1_2_Hs}
\end{align}
\end{theorem}

Using Theorem~\ref{thm:tran_Ak_km1_h}, $\beps_{\B_1}^{(k+1)}$ and $\beps_{\B_2}^{(k+1)}$ in \eqref{eq:eps_1_2} can be represented as
\begin{subequations} \label{eq:approx_eps12}
\begin{align}
\bar{\beps}_{\B_1}^{(k+1)} &= -\Gamma_2^{(k)} A_2^{(k-1)} \beps_{\B_2}^{(k)}  + \Gamma_2^{(k)}  (A_2^{(k-2)} - A_2^{(k-1)}) \h_{\B_2}^{(k-1)}  \nonumber \\
&= -\Gamma_2^{(k)} A_2^{(k-1)} \beps_{\B_2}^{(k)} 
- \Gamma_2^{(k)} P_2 \left( L_2^{(k-1)} \right)^{-1} ( H_{\B_2}^{(k-1)} \beps_{\B_2}^{(k-1)} + H_{\I_2}^{(k-1)} \beps_{\I_2}^{(k-1)} ) \nonumber  \\
&\equiv T_{13}^{(k)} \beps_{\B_2}^{(k)}  + T_{15}^{(k)} \beps_{\B_2}^{(k-1)} + T_{16}^{(k)} \beps_{\I_2}^{(k-1)}  \label{eq:eps1_k}
\end{align} 
and
\begin{align}
\beps_{\B_2}^{(k+1)} 
&= -\overline{\Gamma}_1^{(k+1)} A_1^{(k)} \bar{\beps}_{\B_1}^{(k+1)} 
- \overline{\Gamma}_1^{(k+1)} P_1 \left( L_1^{(k)} \right)^{-1} \left( \overline{H}_{\B_1}^{(k)} \bar{\beps}_{\B_1}^{(k)} + \overline{H}_{\I_1}^{(k)} \bar{\beps}_{\I_1}^{(k)} \right) \nonumber \\
&\equiv S_{31}^{(k)} \bar{\beps}_{\B_1}^{(k+1)} + T_{31}^{(k)} \bar{\beps}_{\B_1}^{(k)} + T_{32}^{(k)} \bar{\beps}_{\I_1}^{(k)}.
\end{align}
\end{subequations}
On the other hand, from \eqref{24}, we have
\begin{align*}
\beps_{\I_{s}}^{(k+1)}  
&=    
\h_{\I_{s}}^{(k+1)}  
- 
\h_{\I_{s}}^{(k)}    =    
\widehat{A}_s^{(k)} \h_{\B_s}^{(k+1)} -  \widehat{A}_s^{(k-1)} \h_{\B_s}^{(k)}   \nonumber \\
&=   
\widehat{A}_s^{(k)} \beps_{\B_s}^{(k+1)} +  \left( \widehat{A}_s^{(k)} - \widehat{A}_s^{(k-1)} \right)\h_{\B_s}^{(k)} \nonumber \\   
&=   
\widehat{A}_s^{(k)} \beps_{\B_s}^{(k+1)} +  \left( L_s^{(k)} \right)^{-1} ({H}_{\B_s}^{(k)}  \beps_{\B_s}^{(k)} + {H}_{\I_s}^{(k)}  \beps_{\I_s}^{(k)}),
\end{align*}
which implies 
\begin{subequations} \label{eq:eps_Is}
\begin{align}
       \bar{\beps}_{\I_{1}}^{(k+1)} -  \widehat{A}_1^{(k)} \bar{\beps}_{\B_1}^{(k+1)} 
      &= \left( L_1^{(k)} \right)^{-1} ( \overline{H}_{\B_1}^{(k)}  \bar{\beps}_{\B_1}^{(k)} + \overline{H}_{\I_1}^{(k)} \bar{\beps}_{\I_1}^{(k)}) \nonumber \\
      &\equiv T_{21}^{(k)} \bar{\beps}_{\B_1}^{(k)} + T_{22}^{(k)} \bar{\beps}_{\I_1}^{(k)} \label{eq:approx_eps2} 
\end{align}
and
\begin{align}
      \beps_{\I_{2}}^{(k+1)} - \widehat{A}_2^{(k)} \beps_{\B_2}^{(k+1)}  &= \left( L_2^{(k)} \right)^{-1} ({H}_{\B_2}^{(k)}  \beps_{\B_2}^{(k)} + {H}_{\I_2}^{(k)}  \beps_{\I_2}^{(k)}) \nonumber \\
      &\equiv T_{43}^{(k)} \beps_{\B_2}^{(k)} + T_{44}^{(k)} \beps_{\I_2}^{(k)}.
\end{align}
\end{subequations}
Combining the results in \eqref{eq:approx_eps12} and \eqref{eq:eps_Is}, we get following theorem.

\begin{theorem} \label{thm:iterative_error_vector}
Let $S_{31}^{(k)}$ and $T_{ij}^{(k)}$ be defined in \eqref{eq:approx_eps12} and \eqref{eq:eps_Is}. Define
\begin{subequations} \label{eq:iter_mtx_TS}
\begin{align}
      \bfzeta^{(k)} &= \begin{bmatrix}
            \bar{\beps}_{\B_1}^{(k)} \\
            \bar{\beps}_{\I_{1}}^{(k)} \\
            \beps_{\B_2}^{(k)} \\
            \beps_{\I_2}^{(k)} \\
            \beps_{\B_2}^{(k-1)} \\
            \beps_{\I_2}^{(k-1)}
      \end{bmatrix}, \quad S^{(k)} = \begin{bmatrix}
           I_{m_1} & 0 & 0 & 0 & 0 & 0 \\
           -  \widehat{A}_1^{(k)} & I_{n_1} & 0 & 0 & 0 & 0 \\
           - S_{31}^{(k)}  & 0 & I_{m_2} & 0 & 0 & 0 \\
           0 & 0 & - \widehat{A}_2^{(k)} & I_{n_2} & 0 & 0 \\
           0 & 0 & 0 & 0 & I_{m_2} & 0 \\
           0 & 0 & 0 & 0 & 0 & I_{n_2}
      \end{bmatrix}  \in \mathbb{C}^{\ell \times \ell}, \\
      T^{(k)} &= \begin{bmatrix}
           0 & 0 & T_{13}^{(k)} & 0 & T_{15}^{(k)} & T_{16}^{(k)} \\
           T_{21}^{(k)} & T_{22}^{(k)} & 0 & 0 & 0 & 0 \\
           T_{31}^{(k)} & T_{32}^{(k)} & 0 & 0 & 0 & 0 \\
           0 & 0 & T_{43}^{(k)} & T_{44}^{(k)} & 0 & 0  \\
           0 & 0 & I_{m_2} & 0 & 0 & 0 \\
           0 & 0 & 0 & I_{n_2} & 0 & 0
      \end{bmatrix} \in \mathbb{C}^{\ell \times \ell}
\end{align}
\end{subequations}
with $\ell = n_1+m_1+2n_2+2m_2$. Then
\begin{align}
\label{3.21}
      \bfzeta^{(k+1)} = (S^{(k)})^{-1} T^{(k)} \bfzeta^{(k)} \equiv \mathcal{T}^{(k)} \bfzeta^{(k)}.
\end{align}
\end{theorem}

\begin{remark}
According to the numerical validation for the quantity of $\mathcal{T}^{(k)}$ in \eqref{3.21}, in Figure~\ref{fig:max_val_spectral_radius_iter_T} of Section~\ref{sec:numerical_validation} (see below), we observe that the spectral radius of $(\mathcal{T}^{(k)} \cdots \mathcal{T}^{(0)})$ and $\max_{k}\max_{ij}(|(\mathcal{T}^{(k)} \cdots \mathcal{T}^{(0)})_{ij}|)$ are all bounded above, for $k$ sufficiently large. 
Therefore, we fairly assume that the sequence $\left\{\left(\mathcal{T}^{(k)}\cdots \mathcal{T}^{(0)}\right)\right\}_{k=0}^{\infty}$ is uniformly bounded and almost diagonalizable. Let $\left(\mathcal{T}^{(k)}\cdots \mathcal{T}^{(0)}\right)^{1/k}$ be the $k$--root of $\left(\mathcal{T}^{(k)}\cdots \mathcal{T}^{(0)}\right)$ for which the eigenvalues are chosen to be the principal values of the eigenvalues of $\left(\mathcal{T}^{(k)}\cdots \mathcal{T}^{(0)}\right)$. There exists a convergent subsequence such that 
\begin{align}
\label{3.22}
\lim_{k\to\infty}\left(\mathcal{T}^{(k)}\cdots \mathcal{T}^{(0)}\right)^{1/k}=\mathcal{A}    
\end{align}
Here, for convenience, ``$k$" denotes the subindices ``$k_j$".
\end{remark}

Based on the assumption \eqref{3.22}, in the following Theorem~\ref{thm3.6}, we will show the Algorithm~\ref{alg:SEM_fix_bnd_k} for spherical equiareal parameterizations is either R-linearly convergent, or forms an asymptotically quasi-periodic solution.
\begin{theorem}
\label{thm3.6}
Set $\mathcal{M}$ be a closed genus-zero surface with normalized area of one. 
Suppose that $\mathcal{A}$ in \eqref{3.22} has the eigendecomposition $\mathcal{A}=X\Lambda X^{-1}$, where $\Lambda=\textup{diag}(e^{\mathrm{i}\theta_1}, \cdots, e^{\mathrm{i}\theta_r}, \lambda_1, \cdots, \lambda_t, \mu_1, \cdots, \mu_u)$ with $\{\theta_1, \ldots, \theta_r\}$ being all distinct, $|\lambda_i|<1$, $i=1, \ldots, t$, and $|\mu_j|>1$, $j=1, \ldots, u$, as well as, $\{\mathbf{x}_j\}_{j=1}^r$ and $\{\mathbf{y}_j\}_{j=1}^r$ being the normalized right and left eigenvectors corresponding to $e^{\mathrm{i}\theta_j}$, $j=1,\ldots, r$, respectively. Then the sequences $\{\mathbf{h}_{\mathtt{I}_s}^{(k)}\}$ and $\{\mathbf{h}_{\mathtt{B}_s}^{(k)}\}$ in \eqref{24}, for $s=1,2$ satisfy
\begin{enumerate}
    \item[\em (i)] either 
    \begin{align*}
    \left\|\beps_s^{(k)}\right\|_{\infty}^{1/k}=\left\|
    \begin{bmatrix}
    \mathbf{h}_{\mathtt{I}_s}^{(k)}-\mathbf{h}_{\mathtt{I}_s}^{(k-1)}\\
    \mathbf{h}_{\mathtt{B}_s}^{(k)}-\mathbf{h}_{\mathtt{B}_s}^{(k-1)}
    \end{bmatrix}\right\|_{\infty}^{1/k}<1,\quad \mbox{(R-linear convergence)}
    \end{align*}
    \item[\em (ii)] or
    \begin{align*}
    \begin{bmatrix}
    \beps_1^{(k)}\\
    \beps_2^{(k)}
    \end{bmatrix}=
    \begin{bmatrix}
    I_1&0&0&0\\
    0&I_2&0&0
    \end{bmatrix}\sum_{j=1}^rb_j {\mathbf{x}_j}e^{\mathrm{i}\theta_jk}+\mathcal{O}(\bfeta_k),
    \end{align*}
    for some $b_j\geq 0$, $\bfeta_k\to0$ as $k\to\infty$, where $I_s=\begin{bmatrix}
    0&I_{n_s}\\
    I_{m_s}&0
    \end{bmatrix}$, $s=1, 2$, (Asymptotically quasi-periodic for $k\gg 1$).
\end{enumerate}
\end{theorem}
\begin{proof}
From assumption and \eqref{3.22} it holds that
\begin{align}
\label{3.23}
\left(\mathcal{T}^{(k)}\cdots\mathcal{T}^{(0)}\right)=\left(\mathcal{A}+E_k\right)^k,\quad \left\|E_k\right\|_{\infty}\to 0,\quad \mbox{as $k\to\infty$}.
\end{align}
If $u\neq\emptyset$, by the continuity of eigenvalues, for $k$ sufficiently large, $(\mathcal{A}+E_k)$ has at least one eigenvalue with modulus larger than one. From  \eqref{3.21} and \eqref{3.23} we have
\begin{align}
    \label{3.24}
    \bfzeta^{(k+1)}=\left(\mathcal{T}^{(k)}\cdots\mathcal{T}^{(0)}\right)\bfzeta_0=\left(\mathcal{A}+E_k\right)^k\bfzeta_0.
\end{align}
Then $\|\bfzeta^{(k+1)}\|_{\infty}$ becomes unbounded for $k$ sufficiently large which contradicts the boundedness of $\|\bfzeta^{(k+1)}\|_{\infty}$. So $\mathcal{A}$ has no eigenvalue with modulus larger than one.

\begin{enumerate}
    \item[(i)] If $r=\emptyset$ and $u=\emptyset$, then $\rho(\mathcal{A})<1$. Since $\lim_{k\to\infty}\left\|\mathcal{A}^k\right\|^{1/k}_{\infty}=\rho(\mathcal{A})\equiv\rho<1$, there is an $N>0$ such that $\left\|\mathcal{A}^N\right\|^{1/N}_{\infty}\leq \rho+\epsilon<1$. We define
    \begin{align*}
        \left\|\x\right\|_*=\left\|\x\right\|_{\infty}+\frac{\left\|\mathcal{A}\x\right\|_{\infty}}{\rho+\epsilon}+\cdots+\frac{\left\|\mathcal{A}^{N-1}\x\right\|_{\infty}}{(\rho+\epsilon)^{N-1}}
    \end{align*}
    and then it holds that
    \begin{align*}
        \left\|\mathcal{A}\right\|_*
        =&\max_{\x\neq0}\frac{\left\|\mathcal{A}\x\right\|_*}{\left\|\x\right\|_*}\\
        \leq&\max_{\x\neq0}\left(\left\|\mathcal{A}\x\right\|_{\infty}+\frac{\left\|\mathcal{A}^2\x\right\|_{\infty}}{\rho+\epsilon}+\cdots+\frac{\left\|\mathcal{A}^{N-1}\x\right\|_{\infty}}{(\rho+\epsilon)^{N-1}}+(\rho+\epsilon)\left\|\x\right\|_{\infty}\right)\\
        =&\rho+\epsilon<1.
    \end{align*}
    There exists a constant $M_{\infty}>0$ such that
    \begin{align*}
        \left\|\bfzeta^{(k+1)}\right\|_{\infty}
        \leq& M_{\infty}\left\|\left(\mathcal{T}^{(k)}\cdots\mathcal{T}^{(0)}\right)\bfzeta_0\right\|_*\\
        =&M_{\infty}\left\|\left(\mathcal{A}+E_k\right)^k\bfzeta_0\right\|_*\leq M_{\infty}\left\|\left(\mathcal{A}+E_k\right)\right\|_*^k\left\|\bfzeta_0\right\|_*
    \end{align*}
    From \eqref{eq:eta_k} and \eqref{eq:iter_mtx_TS} follows that for $k$ sufficiently large
    \begin{align*}
        \left\|\beps^{(k+1)}_s\right\|^{1/k}_{\infty}
        \leq&\left\|\bfzeta^{(k+1)}\right\|^{1/k}_{\infty}\leq M_{\infty}^{1/k}\left(\left\|\mathcal{A}\right\|_*+\left\|E_k\right\|_*\right)\left\|\bfzeta_0\right\|^{1/k}_*<1.
    \end{align*}
    \item[(ii)] If $r\neq\emptyset$ and $s\neq\emptyset$, then
    \begin{align}\label{3.25}
    \bfzeta^{(k+1)}=\left(\mathcal{T}^{(k)}\cdots\mathcal{T}^{(0)}\right)\bfzeta^{(0)}=\left[X(\epsilon_k){\Lambda(\epsilon_k)}^k{X(\epsilon_k)}^{-1}\right]\bfzeta^{(0)},
    \end{align}
    where $\mathcal{A}+E_k=X(\epsilon_k)\Lambda(\epsilon_k){X(\epsilon_k)}^{-1}$ with $\epsilon_k\to 0$, $X(\epsilon_k)\to X$, $\Lambda(\epsilon_k)\to\Lambda$, as $k\to\infty$.
    
    Rewrite \eqref{3.25} as
    \begin{align}\label{3.26}
        \bfzeta^{(k+1)}=X_r(\epsilon_k)\diag\left(e^{\mathrm{i}\theta_jk}\left(1+\epsilon_{k, j}\right)^k\right)_{j=1}^rY_r^{\top}(\epsilon_k)\bfzeta^{(0)}+\mathcal{O}(\hat{\bfeta}_k),
    \end{align}
    where $X_r(\epsilon_k)\to X_r=[\mathbf{x}_1,\cdots, \mathbf{x}_r]$ and $Y_r(\epsilon_k)\to Y_r=[\mathbf{y}_1,\cdots, \mathbf{y}_r]$.
    
    Because of the boundedness of $\left\|\bfzeta^{(k+1)}\right\|_{\infty}$, $\left(1+\epsilon_{k, j}\right)^k$ can not be divergent and should converge to $\beta_j\geq 0$, as $k\to\infty$. Let
    \begin{align}\label{3.27}
        [b_1,\cdots, b_r]^{\top}=\textup{diag}\left(\beta_1,\cdots, \beta_r\right)Y_r^{\top}\bfzeta_0
    \end{align}
    From \eqref{eq:iter_mtx_TS}, \eqref{3.25}-\eqref{3.27} follows that
    \begin{align*}
    \begin{bmatrix}
    \beps_1^{(k)}\\
    \beps_2^{(k)}
    \end{bmatrix}=
    \begin{bmatrix}
    I_1&0&0&0\\
    0&I_2&0&0
    \end{bmatrix}\sum_{j=1}^rb_j\mathbf{x}_je^{\mathrm{i}\theta_jk}+\mathcal{O}(\bfeta_k),
    \end{align*}
    where $\bfeta_k\to 0$, as $k\to\infty$.
\end{enumerate}
\end{proof}

\begin{remark}
In practice, case (ii) of Theorem~\ref{thm3.6} rarely occurs. We observe that it can happen if the boundary 
$\h_{\B_s}^{(k)}$ is not positioned during the iteration, then $\h_{\B_s}^{(k)}$ and $\h_{\B_s}^{(k-1)}$ may differ by the sum of some vectors $\{ \mathbf{x}_j \}$ with rotations $\{e^{\mathrm{i}\theta_j k}\}$ at each other iteration, which acts asymptomatically as a quasi-periodic solution. 
\end{remark} 

\begin{figure}
\centering
\begin{tabular}{cccc}
David & Venus & apple & brain \\
\includegraphics[height=3cm]{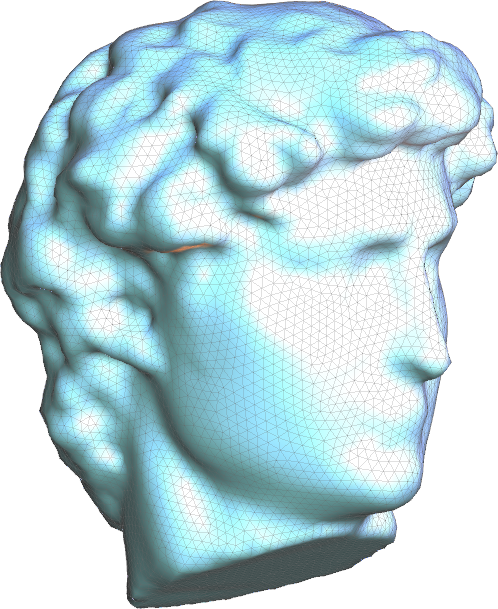} &
\includegraphics[height=3cm]{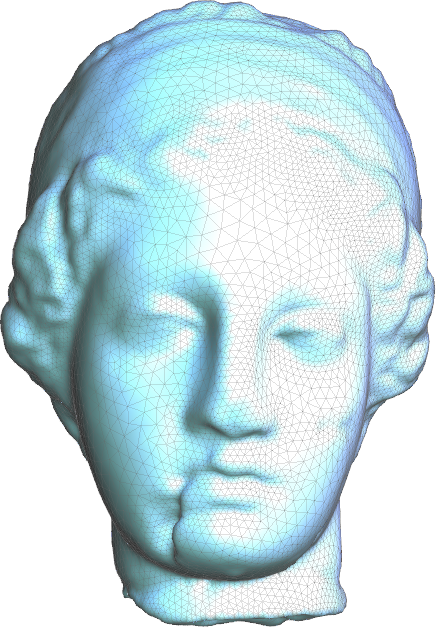} &
\includegraphics[height=3cm]{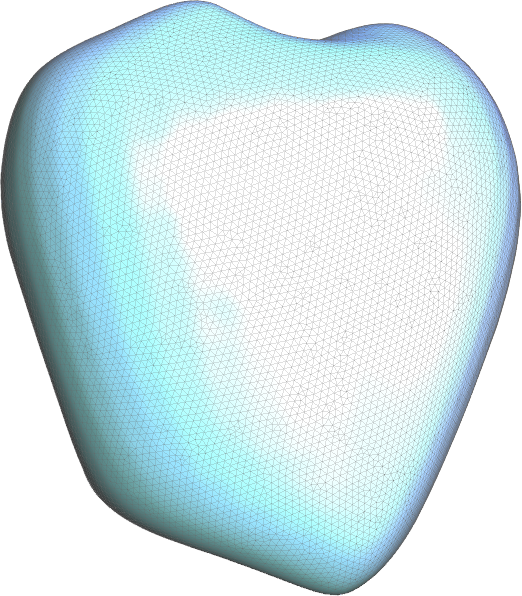} & 
\includegraphics[height=3cm]{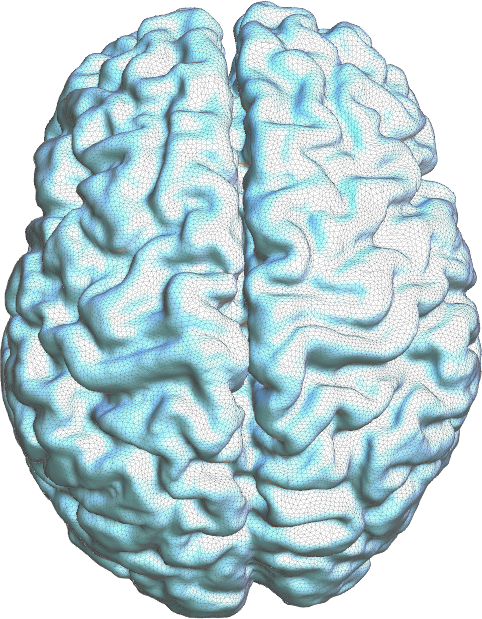}
\end{tabular}
\caption{The benchmark genus-zero closed surface models.}
\label{fig:benchmark_mesh}
\end{figure}

\section{Numerical validations and experiments} \label{sec:4}
In this section, we first give the numerical study on the boundness of matrices $\{(\mathcal{T}^{(k)}\cdots \mathcal{T}^{(0)})\}$ for $k$ large in Subsection~\ref{sec:numerical_validation}. Then, to check the validity of Theorem~\ref{thm3.6}, in Subsection~\ref{sec:numerical_results} we demonstrate the numerical experiment on the R-linear convergence  for various closed genus-zero surface models as in Figure~\ref{fig:benchmark_mesh} from GitHub~\cite{GitHub}, Gu's website~\cite{GuWeb}, and the BraTS 2021 databases~\cite{BaGB21,BaAS17}. The associated matrix sizes $n=\#\V(\M)$,  $n_s = \# \mathtt{I}_{s}$ and $m_s = \# \mathtt{B}_s$, for $s=1, 2$, as in \eqref{eq:idx_sets_I_Bprime_B} with three different mesh sizes are displayed in Table~\ref{tab:matrix_size}. The maximal length of the meshes are $0.02$, $0.015$ and $0.01$ for ``mesh1'', ``mesh2'' and ``mesh3'', respectively.  

\begin{table}[h]
  \centering
  \begin{tabular}{|c|c|r|r|c|r|c|c|} \hline
     \multicolumn{2}{ |c| }{}& $n$ & $n_1$ & $m_1$ & $n_2$ & $m_2$ \\ \hline 
     \multirow{3}{*}{apple} & mesh1 & 10316 & 5207 & 182 & 5201 & 183 \\ 
     & mesh2 & 18269 & 9231 & 251 & 9228 & 251 \\ 
     & mesh3  & 40777 & 20632 & 398 & 20554 & 397 \\ \hline 
     \multirow{3}{*}{brain} & mesh1 & 9699 & 4974 & 182 & 4967 & 184 \\ 
     & mesh2 & 17253 & 8791 & 256 & 8896 & 250 \\ 
     & mesh3  & 38587 & 19736 & 356 & 19793 & 354 \\ \hline 
     \multirow{3}{*}{David} & mesh1 & 9691 & 4895 & 169 & 4889 & 170 \\ 
     & mesh2 & 17125 & 8741 & 230 & 8714 & 224 \\ 
     & mesh3  & 38484 & 19636 & 337 & 19598 & 339 \\ \hline 
     \multirow{3}{*}{Venus} & mesh1 & 8768 & 4420 & 183 & 4429 & 185 \\ 
     & mesh2 & 15543 & 7876 & 219 & 7816 & 218 \\ 
     & mesh3  & 34807 & 17629 & 378 & 17521 & 370  \\ \hline 
  \end{tabular} 
  \caption{Dimensions with various mesh sizes of $n$, $n_s$ and $m_s$, for $s = 1, 2$.} \label{tab:matrix_size}
\end{table}

\subsection{Numerical validations} \label{sec:numerical_validation}
\begin{figure}[tbhp]
\center
%
%
\begin{subfigure}[b]{0.49\textwidth}
\center
    \includegraphics[width=\textwidth]{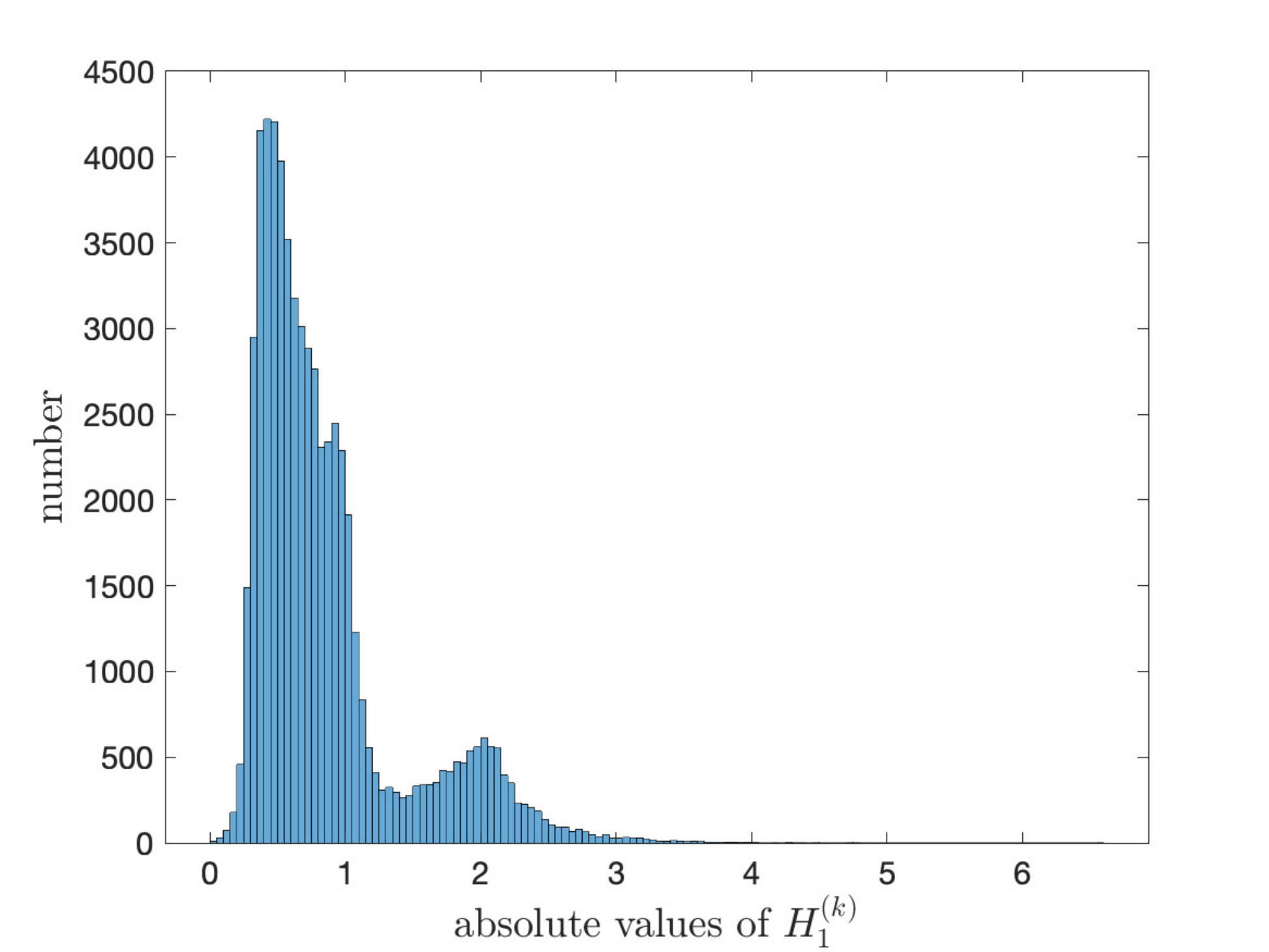}
\caption{$|H_1^{(k)}|$}
    \label{fig:hist_H1_k_brain1}
\end{subfigure}
\begin{subfigure}[b]{0.49\textwidth}
\center
    \includegraphics[width=\textwidth]{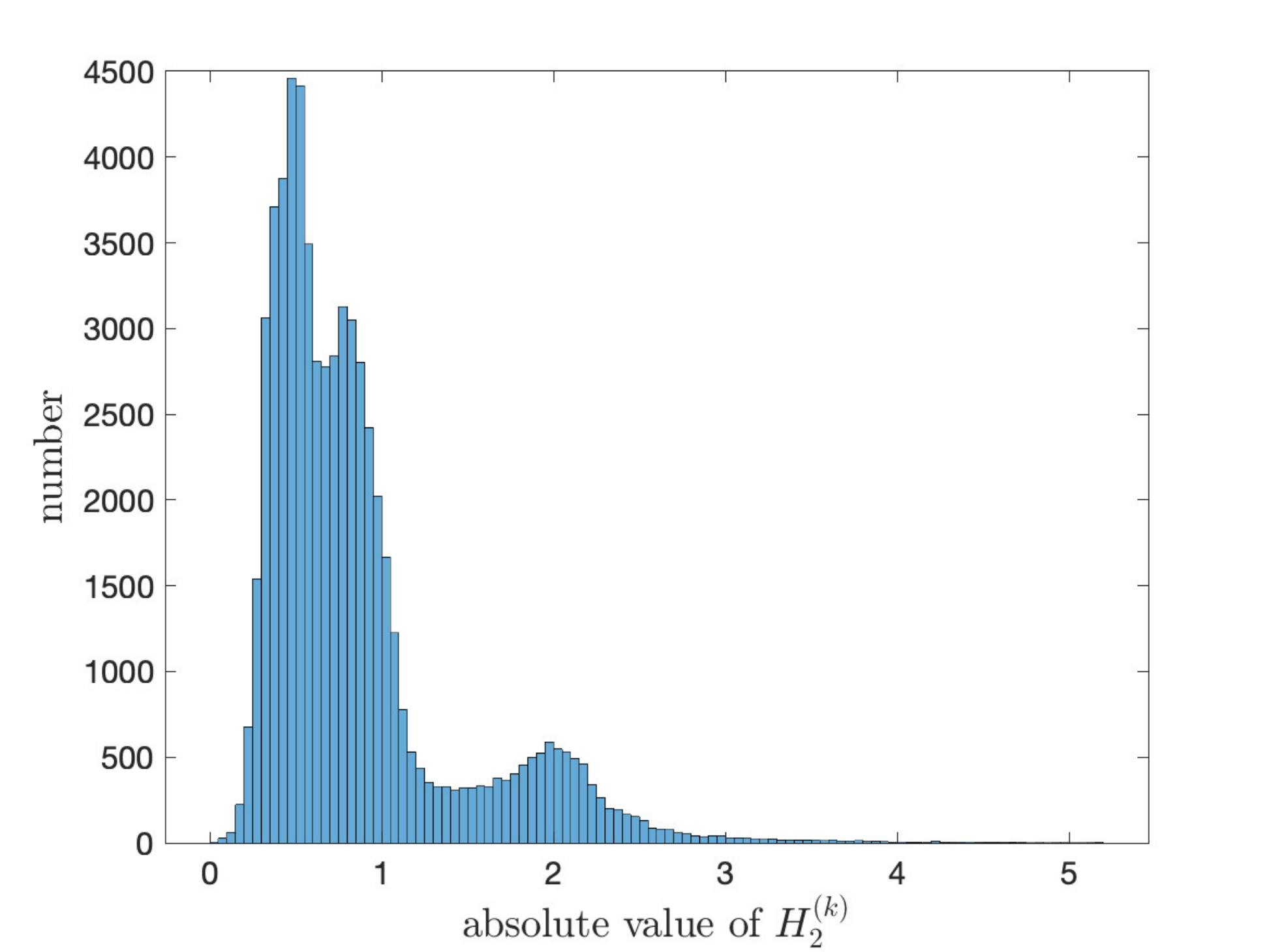}
\caption{$|H_2^{(k)}|$}
    \label{fig:hist_H2_km1_brain1}
\end{subfigure}
\caption{Histograms of the nonzero  elements of $|H_1^{(k)}|$ and $| H_2^{(k)}|$ in (A) and (B), respectively.}
\label{fig:hist_H12_L12_brain1}
\end{figure}

In this subsection, the numerical results are produced by using mesh size ``mesh3'' in Table~\ref{tab:matrix_size}.
To verify the boundedness of the matrices $\{(\mathcal{T}^{(k)}\cdots \mathcal{T}^{(0)})\}$ in Theorem~\ref{thm:tran_Ak_km1_h},  we first construct two crucial matrices $H_{\I_s}^{(k)}$ and  $H_{\B_s}^{(k)}$ from a brain image model of \cite{BaGB21,BaAS17} such that the equality in \eqref{eq:L_k-L_km1_2_Hs} holds. The magnitudes of $\widehat{w}_{ij}^{(k)}\equiv ( [L_s^{(k)}, B_s^{(k)}] - [L_s^{(k-1)}, B_s^{(k-1)} ])_{ij}$,  $\| \mathbf{g}_s^{(k)}\|$ and $\| \beps_s^{(k)} \|$ can be checked to be $\widehat{w}_{ij}^{(k)} = \mathcal{O}(\delta)$ (small), $\| \mathbf{g}_s^{(k)} \| = \mathcal{O}(1)$ and $\| \beps_{s}^{(k)} \| = \mathcal{O}(\delta)$ (small). In Figures~\ref{fig:hist_H1_k_brain1} and \ref{fig:hist_H2_km1_brain1}, we demonstrate the histograms of the nonzero elements of $H_s^{(k)} \equiv [H_{\I_s}^{(k)}, H_{\B_s}^{(k)} ]$ for $k = 4$, $s=1, 2$, respectively. We see that the maximal absolute elements are all $\mathcal{O}(1)$ and the most of the absolute elements are in the interval $[0.3, 1]$. This means that we have constructed reasonable matrices $H_{\I_s}^{(k)}$ and  $H_{\B_s}^{(k)}$ to satisfy the representation in  \eqref{eq:L_k-L_km1_2_Hs}.

\begin{figure}
\center
\begin{subfigure}[b]{0.32\textwidth}
\center
    \includegraphics[width=\textwidth]{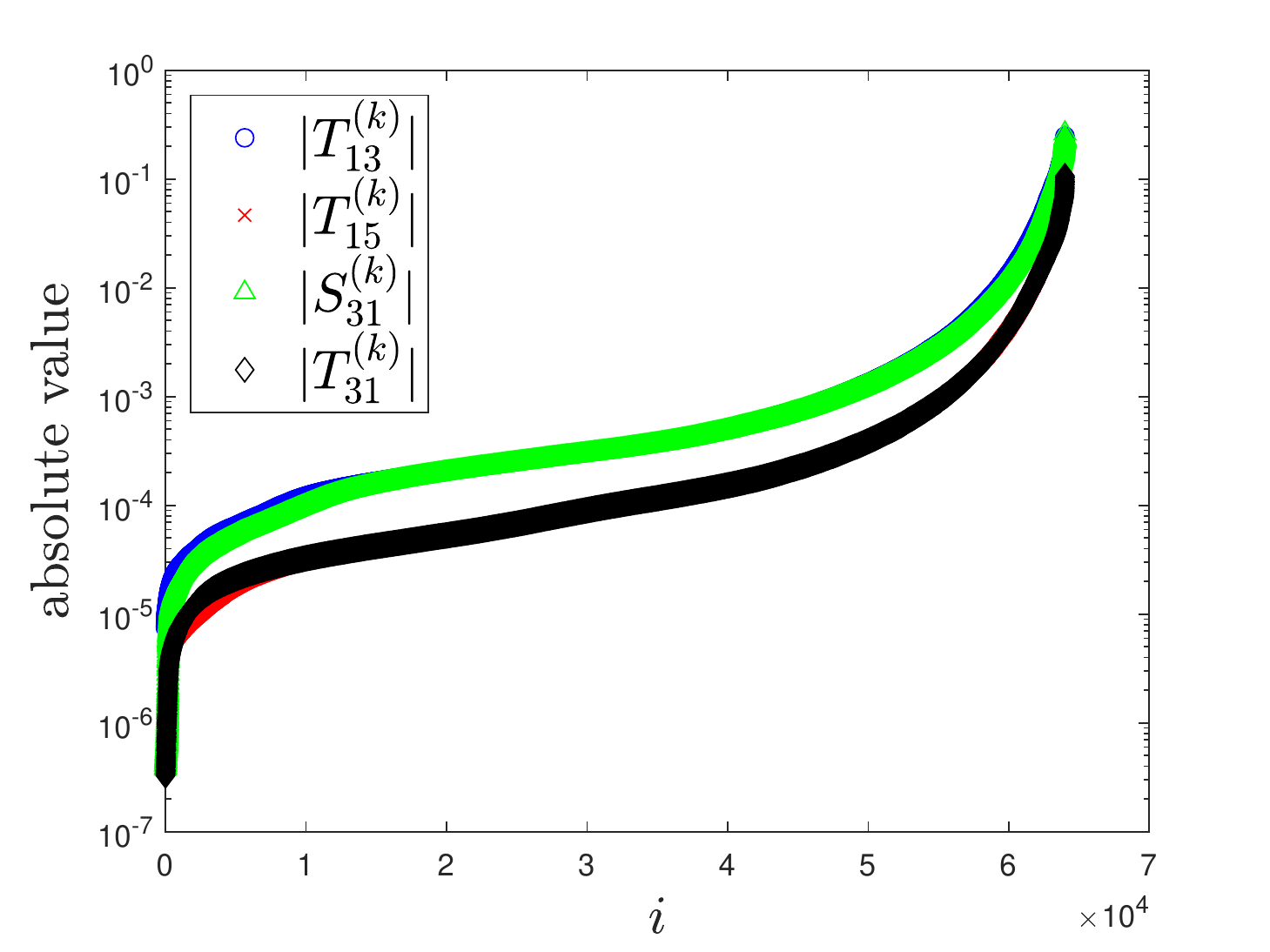}
\caption{$|T_{13}^{(k)}|$, $|T_{15}^{(k)}|$, $|S_{31}^{(k)}|$, $|T_{31}^{(k)}|$}
    \label{fig:abs_value_dist_T13_S31}
\end{subfigure}
\begin{subfigure}[b]{0.32\textwidth}
\center
    \includegraphics[width=\textwidth]{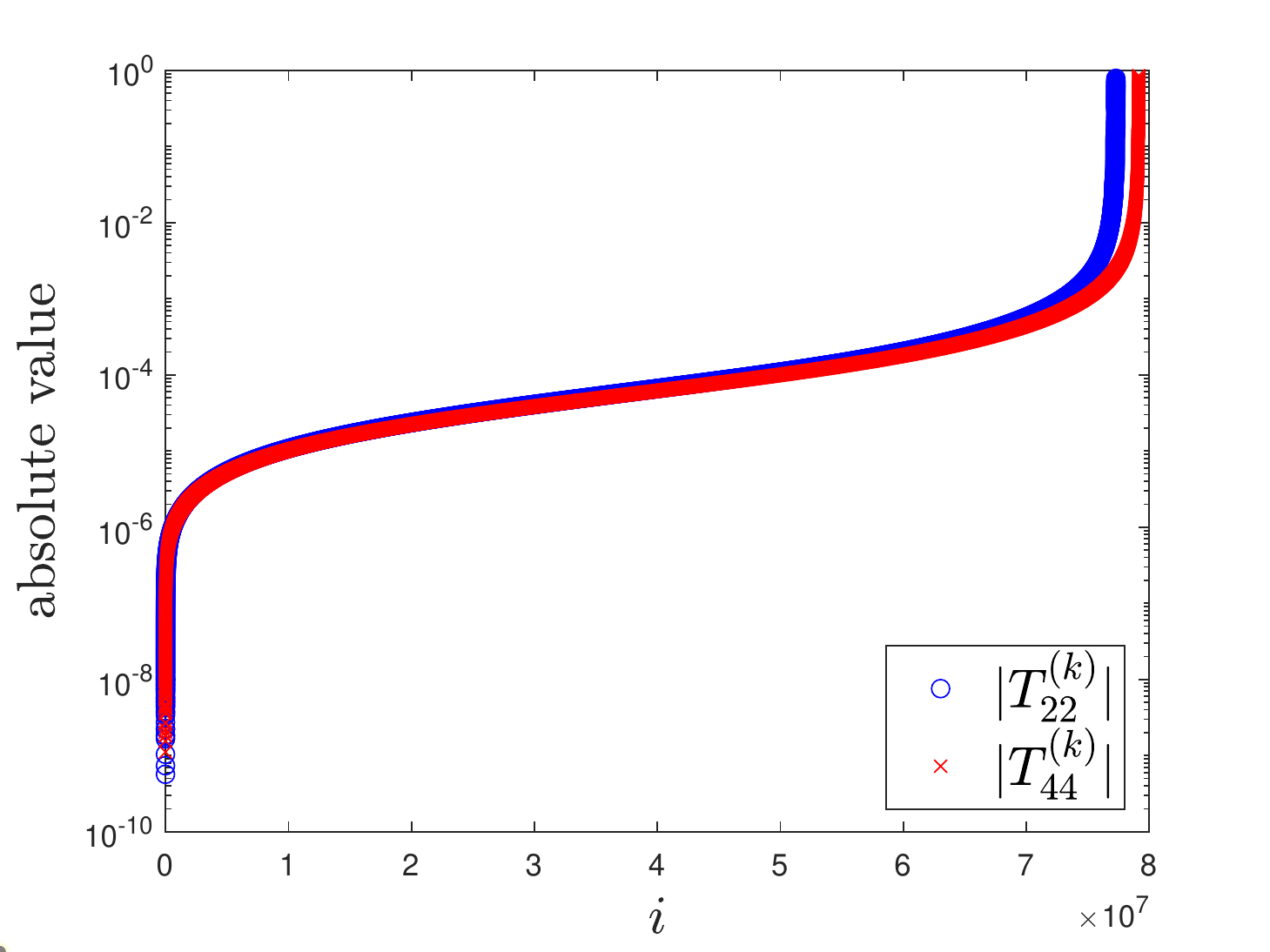}
\caption{$|T_{22}^{(k)}|$, $|T_{44}^{(k)}|$}
    \label{fig:abs_value_dist_T22_T44}
\end{subfigure}
\begin{subfigure}[b]{0.32\textwidth}
\center
    \includegraphics[width=\textwidth]{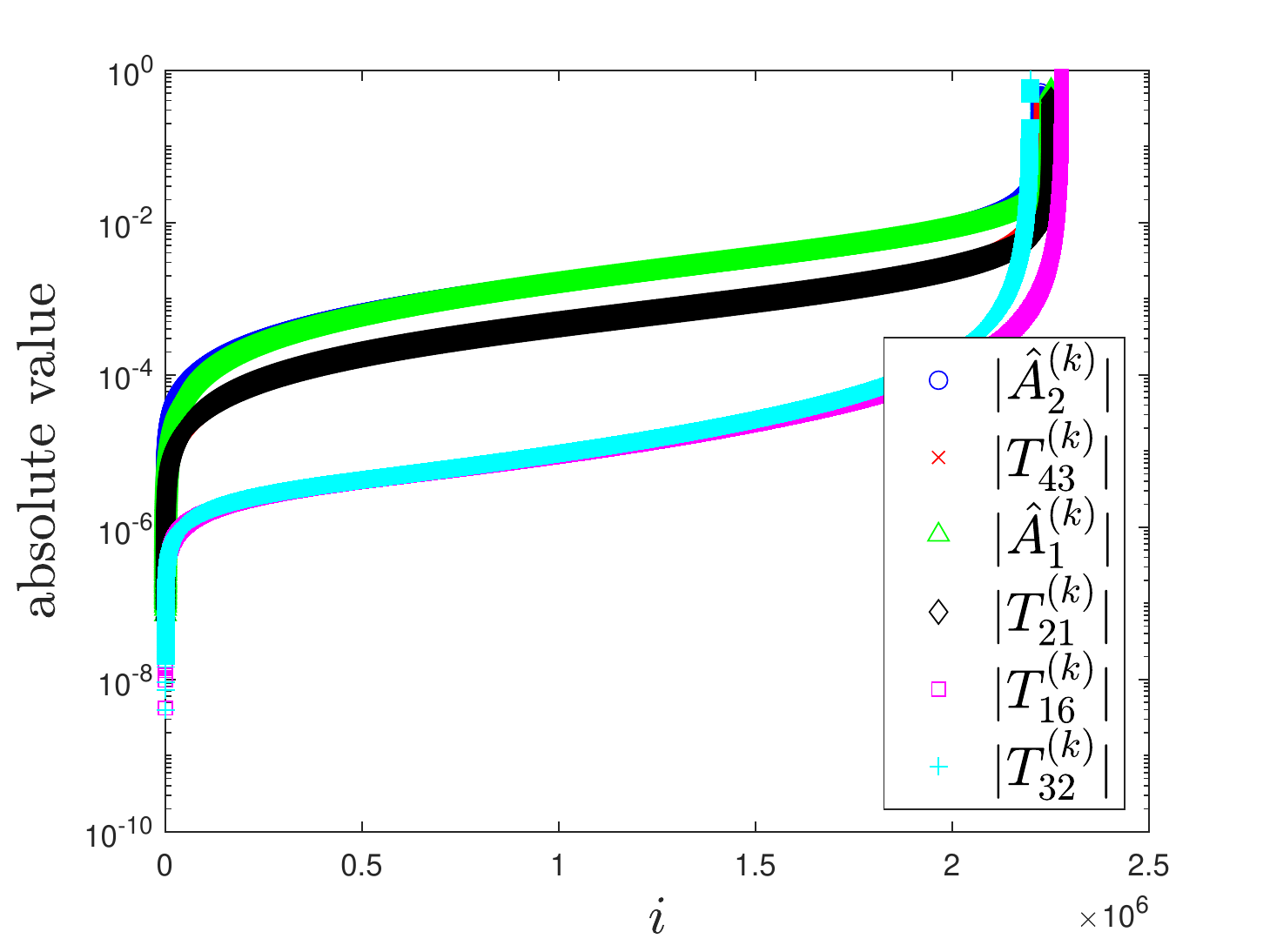}
\caption{$|\widehat{A}_s^{(k)}|$, $|T_{16}^{(k)}|$, $|T_{21}^{(k)}|$, $|T_{32}^{(k)}|$, $|T_{43}^{(k)}|$}
    \label{fig:abs_value_dist_T43_A2}
\end{subfigure} 
\caption{The nonzero  elements (sorting in ascending order) of $|\widehat{A}_s^{(k)}|$, $|S_{31}^{(k)}|$, and $|T_{ij}^{(k)}|$ in \eqref{eq:iter_mtx_TS}.}
\label{fig:abs_value_dist_Tij}
\end{figure}

Next, we plot the distributions of the elements of $|\widehat{A}_s^{(k)}|$, $|S_{31}^{(k)}|$, and $|T_{ij}^{(k)}|$ for $k = 4$ in \eqref{eq:approx_eps12} and \eqref{eq:eps_Is}. We sort the nonzero absolute elements of those matrices in ascending order and display the sorting results in Figure~\ref{fig:abs_value_dist_Tij}. The results show that all the nonzero absolute elements of those matrices are less than 1. Most of the nonzero absolute elements of $\{ \widehat{A}_s^{(k)}, T_{43}^{(k)}, T_{21}^{(k)}\}$, $\{ T_{13}^{(k)}, T_{15}^{(k)}, T_{31}^{(k)}, S_{31}^{(k)} \}$ and $\{ T_{16}^{(k)}, T_{22}^{(k)}, T_{32}^{(k)}, T_{44}^{(k)} \}$ are belonging to the interval $[10^{-4}, 10^{-2}]$, $[10^{-5}, 10^{-3}]$ and $[ 10^{-6}, 10^{-3}]$, respectively. 

\begin{figure}[tbhp]
\center
\begin{subfigure}[b]{0.49\textwidth}
\center
    \includegraphics[width=\textwidth]{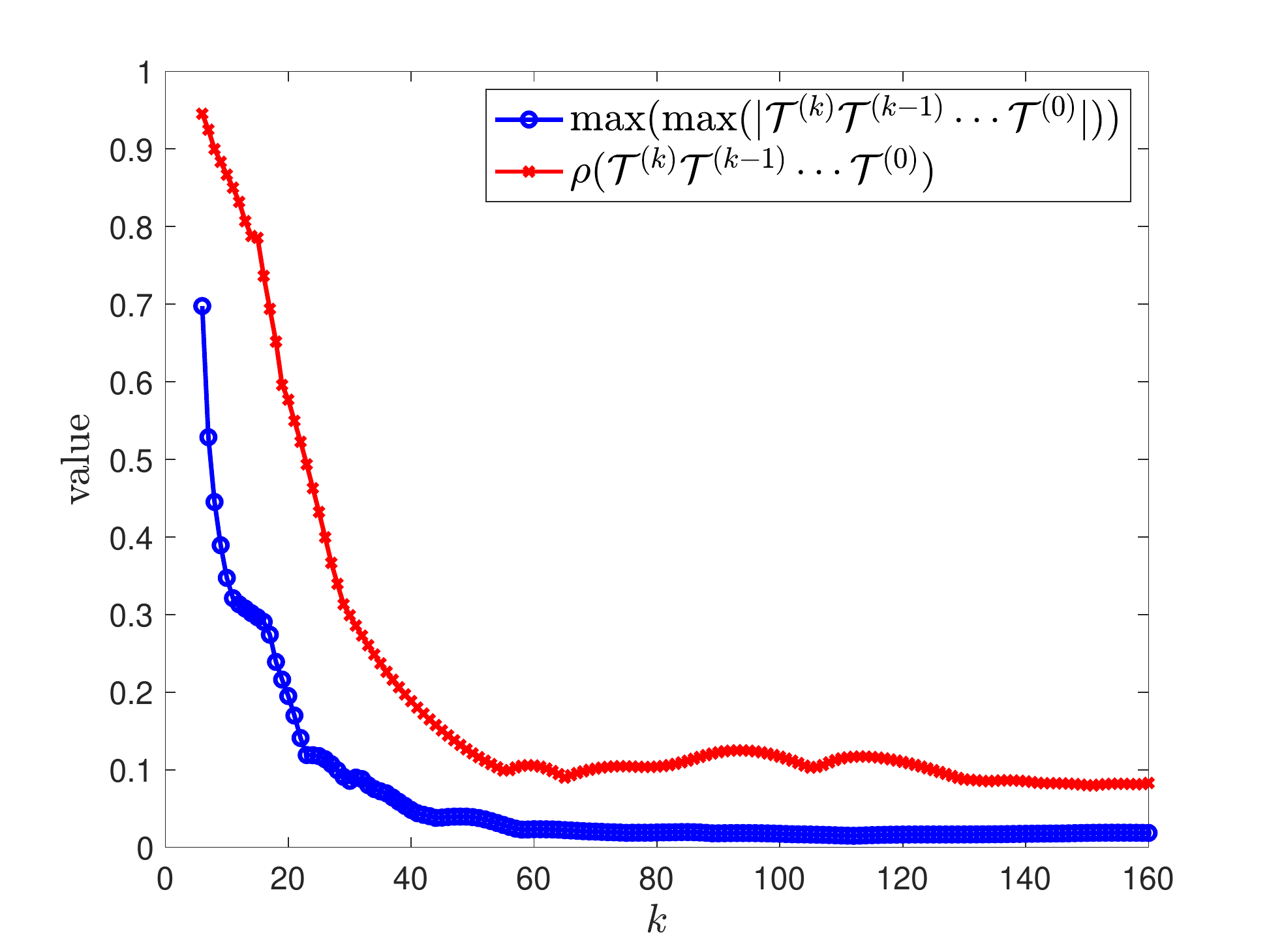}
\caption{brain}
    \label{fig:max_val_spectral_radius_iter_T_brain}
\end{subfigure}
\begin{subfigure}[b]{0.49\textwidth}
\center
    \includegraphics[width=\textwidth]{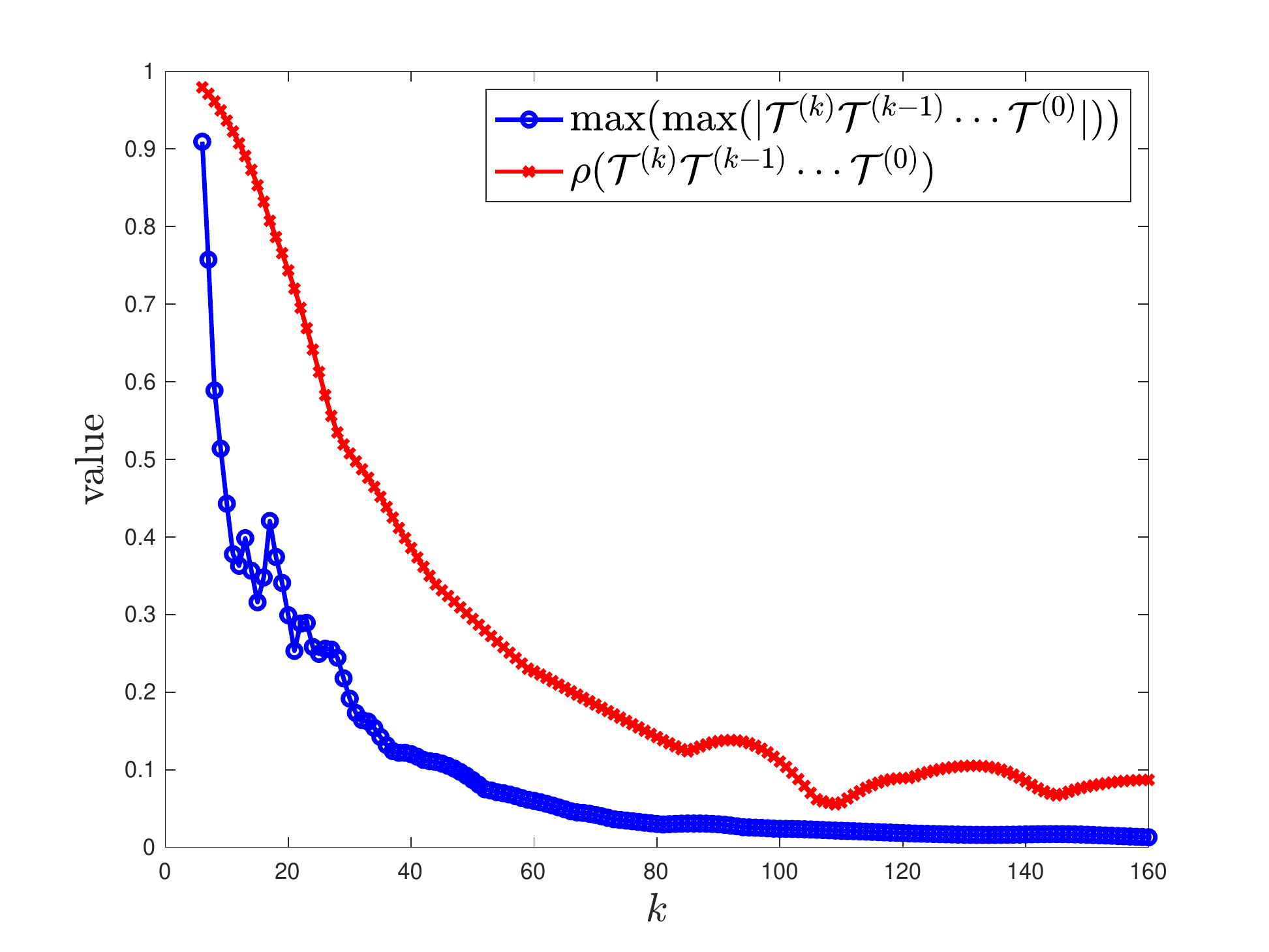}
\caption{Venus}
    \label{fig:max_val_spectral_radius_iter_T_Venus}
\end{subfigure}
\caption{Maximal absolute elements and the spectral radius of the matrix $(\mathcal{T}^{(k)}\cdots \mathcal{T}^{(0)})$ with various $k$ for the brain image and Venus benchmarks.}
\label{fig:max_val_spectral_radius_iter_T}
\end{figure}

The iterative matrix $\mathcal{T}^{(k)}$ in \eqref{3.21} is generated by the matrices $|\widehat{A}_s^{(k)}|$, $|S_{31}^{(k)}|$, and $|T_{ij}^{(k)}|$. In Figure~\ref{fig:max_val_spectral_radius_iter_T}, we demonstrate the maximal absolute elements and the spectral radius of the matrix $(\mathcal{T}^{(k)}\cdots \mathcal{T}^{(0)})$ with various $k$ for the brain image and Venus from benchmarks. Numerical results show that the maximal absolute elements and the spectral radius both are almost decreasing and uniformly bounded for $k$ sufficiently large.  

\subsection{Stretch factors and stretch energies vs mesh sizes} 
\label{sec:numerical_results}
\begin{figure}
\center
\begin{subfigure}[b]{0.49\textwidth}
\center
    \includegraphics[width=\textwidth]{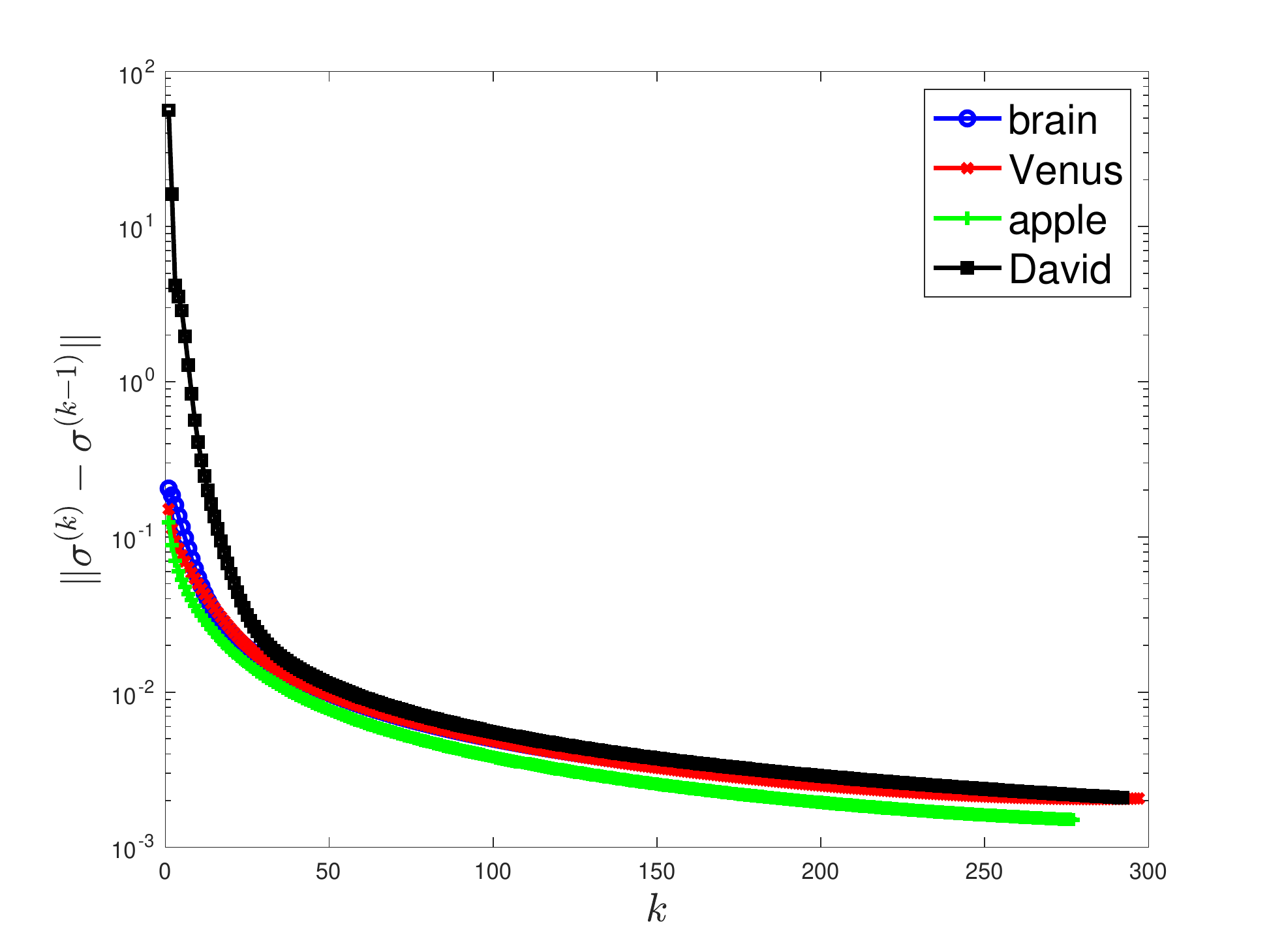}
\caption{$\| \bsigma^{(k)} - \bsigma^{(k-1)} \|_2$}
    \label{fig:stretch_factor_conv}
\end{subfigure} 
\begin{subfigure}[b]{0.49\textwidth}
\center
    \includegraphics[width=\textwidth]{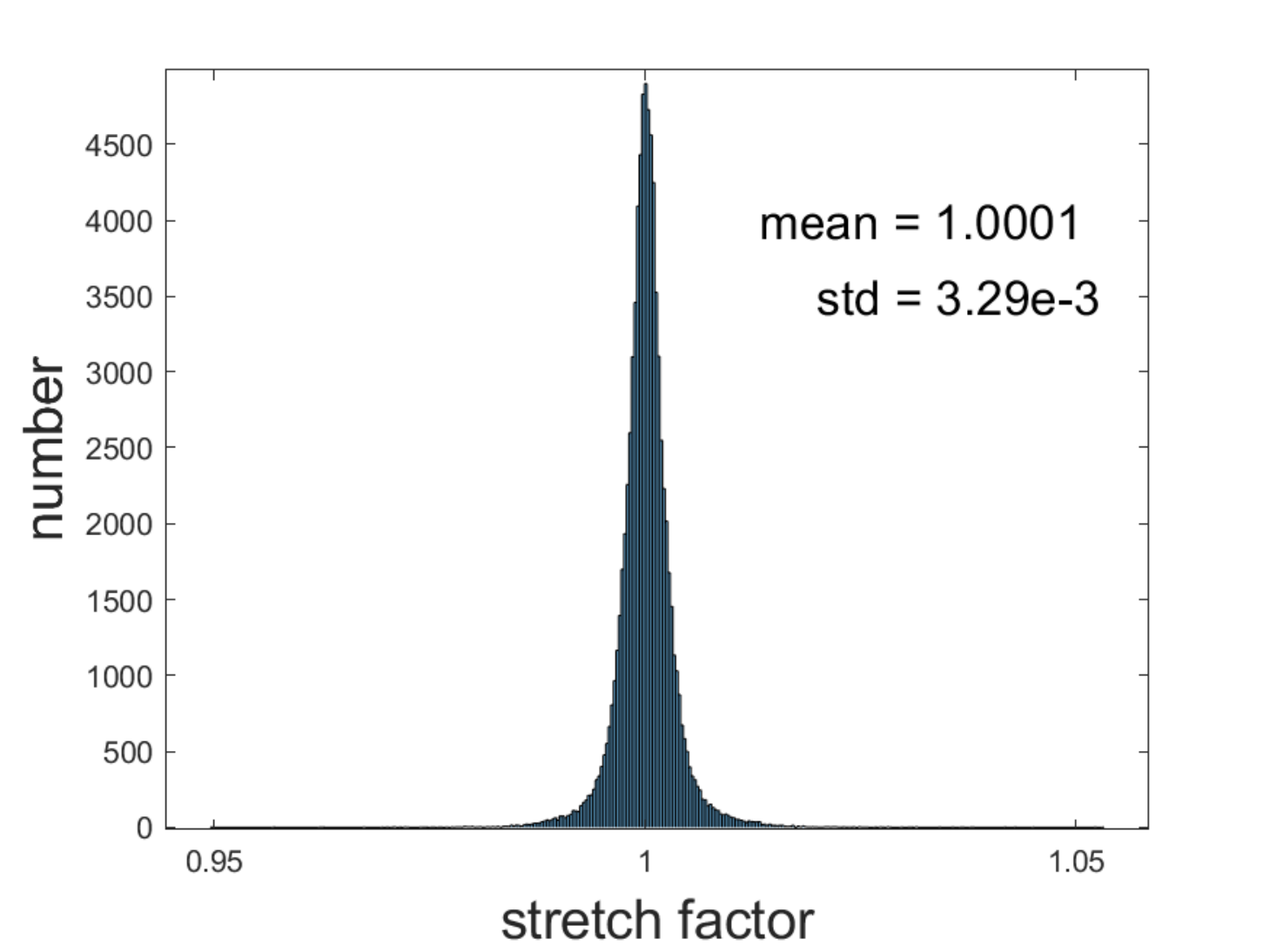}
\caption{apple}
    \label{fig:stretch_factor_Apple}
\end{subfigure} 
\begin{subfigure}[b]{0.32\textwidth}
\center
    \includegraphics[width=\textwidth]{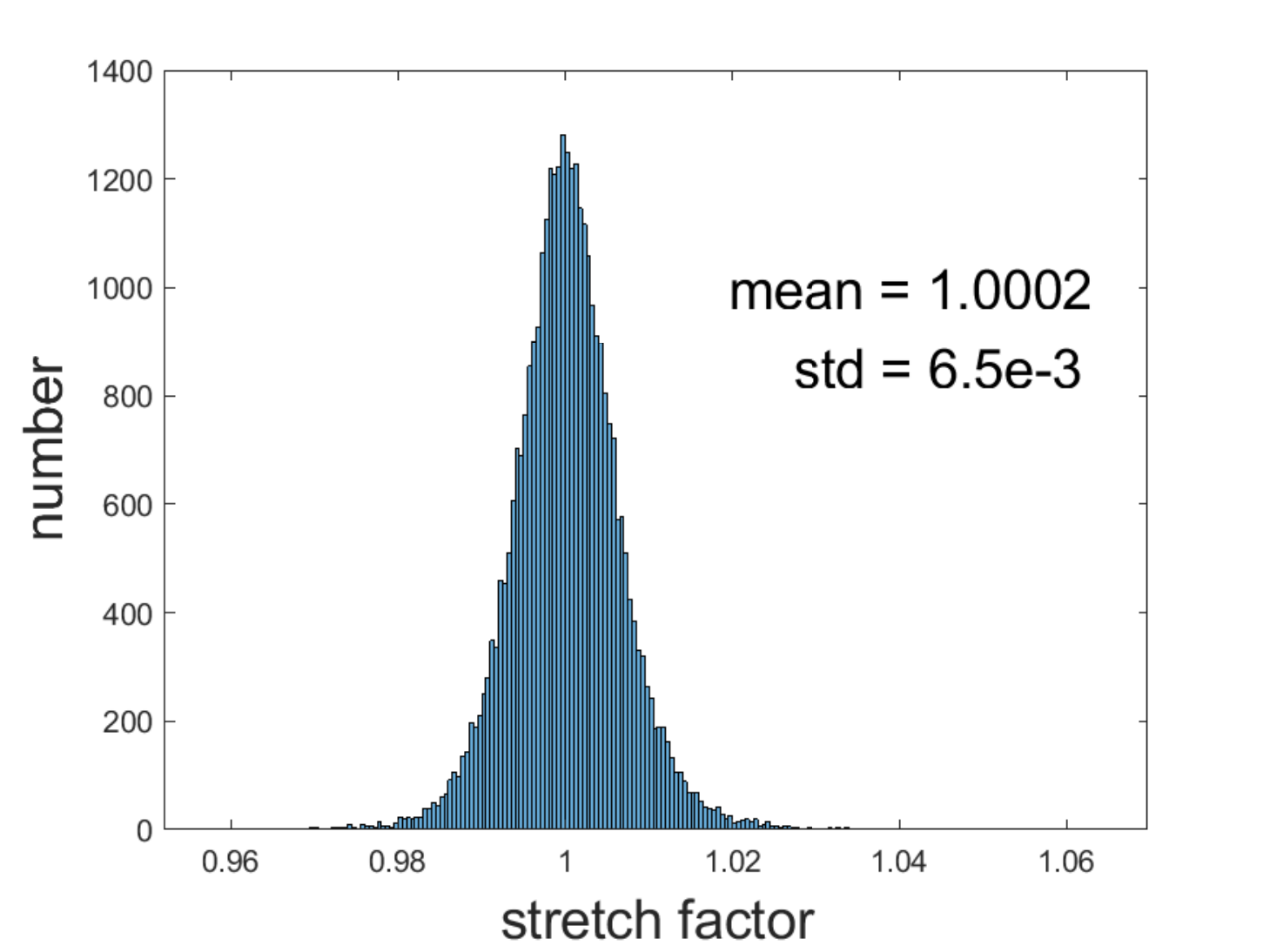}
\caption{brain}
    \label{fig:stretch_factor_brain}
\end{subfigure}
\begin{subfigure}[b]{0.32\textwidth}
\center
    \includegraphics[width=\textwidth]{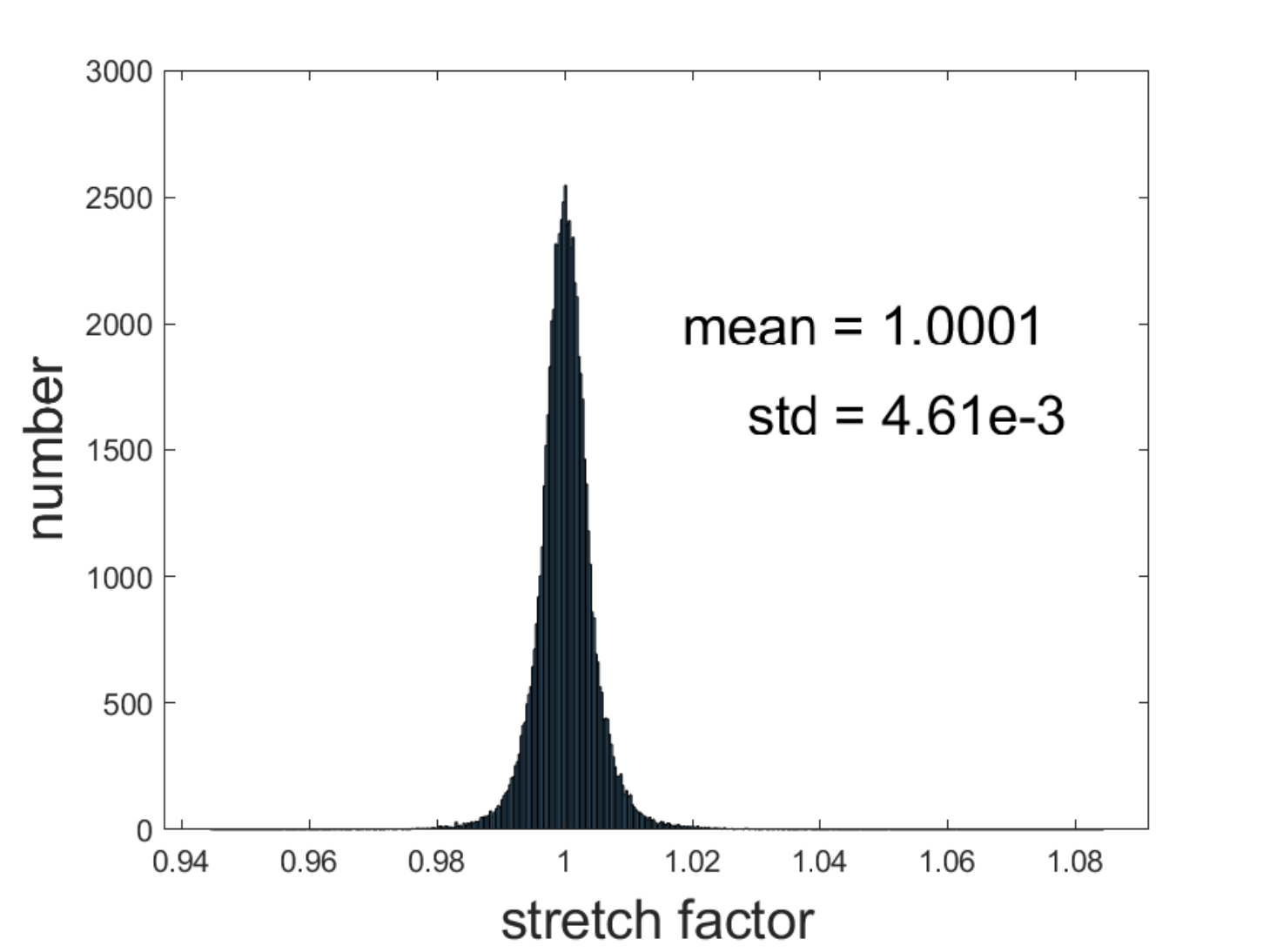}
\caption{Venus}
    \label{fig:stretch_factor_Venus}
\end{subfigure}
\begin{subfigure}[b]{0.32\textwidth}
\center
    \includegraphics[width=\textwidth]{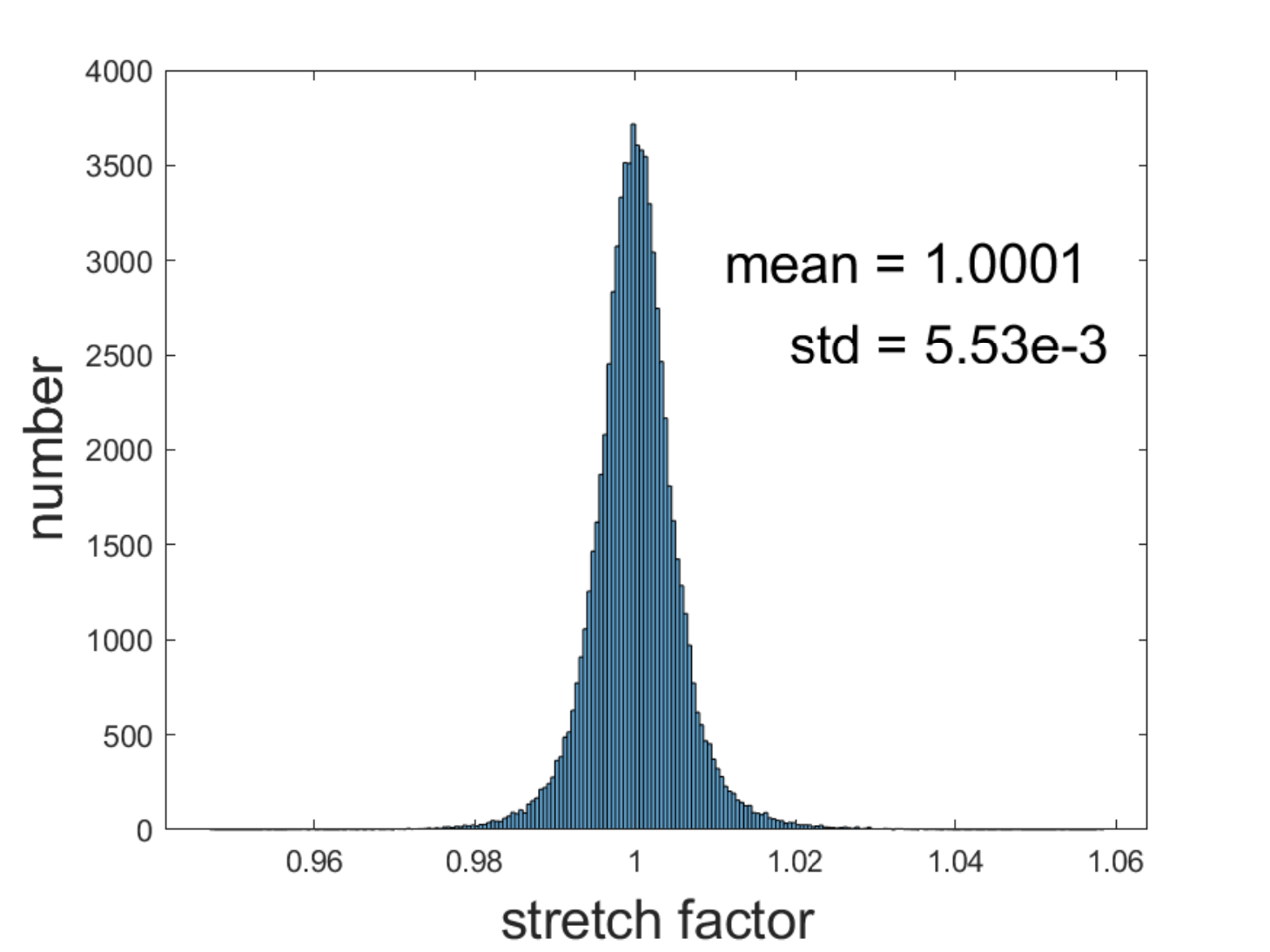}
\caption{David}
    \label{fig:stretch_factor_DavidHead}
\end{subfigure}
\caption{Convergence behavior of $\bsigma^{(k)}$ and histograms of the stretch factors for the four models.}
\label{fig:stretch_factors}
\end{figure}

For the continuous version, Theorem~\ref{Equiareal} shows that the map $\h$ of \eqref{invh} satisfies the equiareal (area-preserving) property.  For the discrete version, we present the area-preserving property for Algorithm~\ref{alg:SEM_fix_bnd_k} by the vector $\bsigma^{(k)}$ of the stretch factor as
\begin{align*}
    \bsigma^{(k)} = \left[ |\Delta_l |/|f^{(k)}(\Delta_l)| \right].
\end{align*}
In Figure \ref{fig:stretch_factor_conv}, we plot the results of $\| \bsigma^{(k)} - \bsigma^{(k-1)} \|_2$ for each iteration $k$ to illustrate the convergence behavior of four benchmark models in Figure~\ref{fig:benchmark_mesh} with mesh size ``mesh3''. When the iteration is convergent, we demonstrate the histograms of the convergent $\bsigma^{(k)}$, the mean value (denoted as ``mean'') and the standard deviation (denoted as ``std'') of $\bsigma^{(k)}$ for the benchmark models in Figures~\ref{fig:stretch_factor_Apple} -- \ref{fig:stretch_factor_DavidHead}. These mean values and  standard deviations show that  stretch factors (local area ratios) are closed to one. This indicates that the map produced by Algorithm~\ref{alg:SEM_fix_bnd_k} preserves the local area ratios very well.

\begin{table}[h]
  \centering
  \begin{tabular}{|c|c|c|c|c|c|c|c|}  \hline
     & \multicolumn{3}{ c| }{apple} & \multicolumn{3}{ c| }{brain}  \\ \cline{2-7}
     & $\mbox{mean} - 1$ & std & $4 \pi - E_{S}(f)$ & $\mbox{mean} - 1$ & std & $4 \pi - E_{S}(f)$  \\ \hline
     mesh1 & $3.4\times 10^{-4}$ & $5.2 \times 10^{-3}$ & $7.5 \times 10^{-3}$ & $4.3 \times 10^{-4}$ & $8.2 \times 10^{-3}$ & $7.4 \times 10^{-3}$  \\
     mesh2 & $2.0\times 10^{-4}$ & $4.9 \times 10^{-3}$ & $4.1 \times 10^{-3}$ & $2.5 \times 10^{-4}$ & $6.6 \times 10^{-3}$ & $4.1 \times 10^{-3}$  \\
     mesh3  & $9.2 \times 10^{-5}$ & $3.3 \times 10^{-3}$ & $1.8 \times 10^{-3}$ & $1.4 \times 10^{-4}$ & $6.5 \times 10^{-3}$ & $1.6 \times 10^{-3}$   \\ \hline \hline
     & \multicolumn{3}{ c| }{David} & \multicolumn{3}{ c| }{Venus}  \\ \cline{2-7}
     & $\mbox{mean} - 1$ & std & $4 \pi - E_{S}(f)$ & $\mbox{mean} - 1$ & std & $4 \pi - E_{S}(f)$  \\ \hline
     mesh1 & $4.9\times 10^{-4}$ & $1.1 \times 10^{-2}$ & $6.8 \times 10^{-3}$ & 
         $4.3 \times 10^{-4}$ & $8.6 \times 10^{-3}$ & $8.2 \times 10^{-3}$  \\
     mesh2 & $2.6\times 10^{-4}$ & $8.0 \times 10^{-3}$ & $3.9 \times 10^{-3}$ & 
         $2.6 \times 10^{-4}$ & $7.4 \times 10^{-3}$ & $4.5 \times 10^{-3}$  \\
     mesh3  & $1.1 \times 10^{-4}$ & $5.5 \times 10^{-3}$ & $1.7 \times 10^{-3}$ & 
         $1.1 \times 10^{-4}$ & $4.6 \times 10^{-3}$ & $2.1 \times 10^{-3}$ \\ \hline
  \end{tabular} 
  \caption{Mean, std and stretch energy $E_S(f)$ in \eqref{17a-1} with various mesh sizes.} \label{tab:conv_energy}
\end{table}

Next, we list the convergence of the stretch factor and stretch energy $E_S(f)$  in \eqref{17a-1} with mesh sizes ``mesh1'', ``mesh2'' and ``mesh3'' in Table~\ref{tab:conv_energy}. The results show that the stretch factor is approaching to $1$ as the mesh size becomes smaller. The results also demonstrate that the stretch energy is increasing and convergent to $4\pi$ (the area of the 2-sphere $\mathbb{S}^2$) for each benchmark problem.  This means that the fine mesh can improve the area-preserving property and the accuracy of the stretch energy.

\subsection{Asymptotically R-linear convergence}
\begin{figure}
\center
\begin{subfigure}[b]{0.49\textwidth}
\center
    \includegraphics[width=\textwidth]{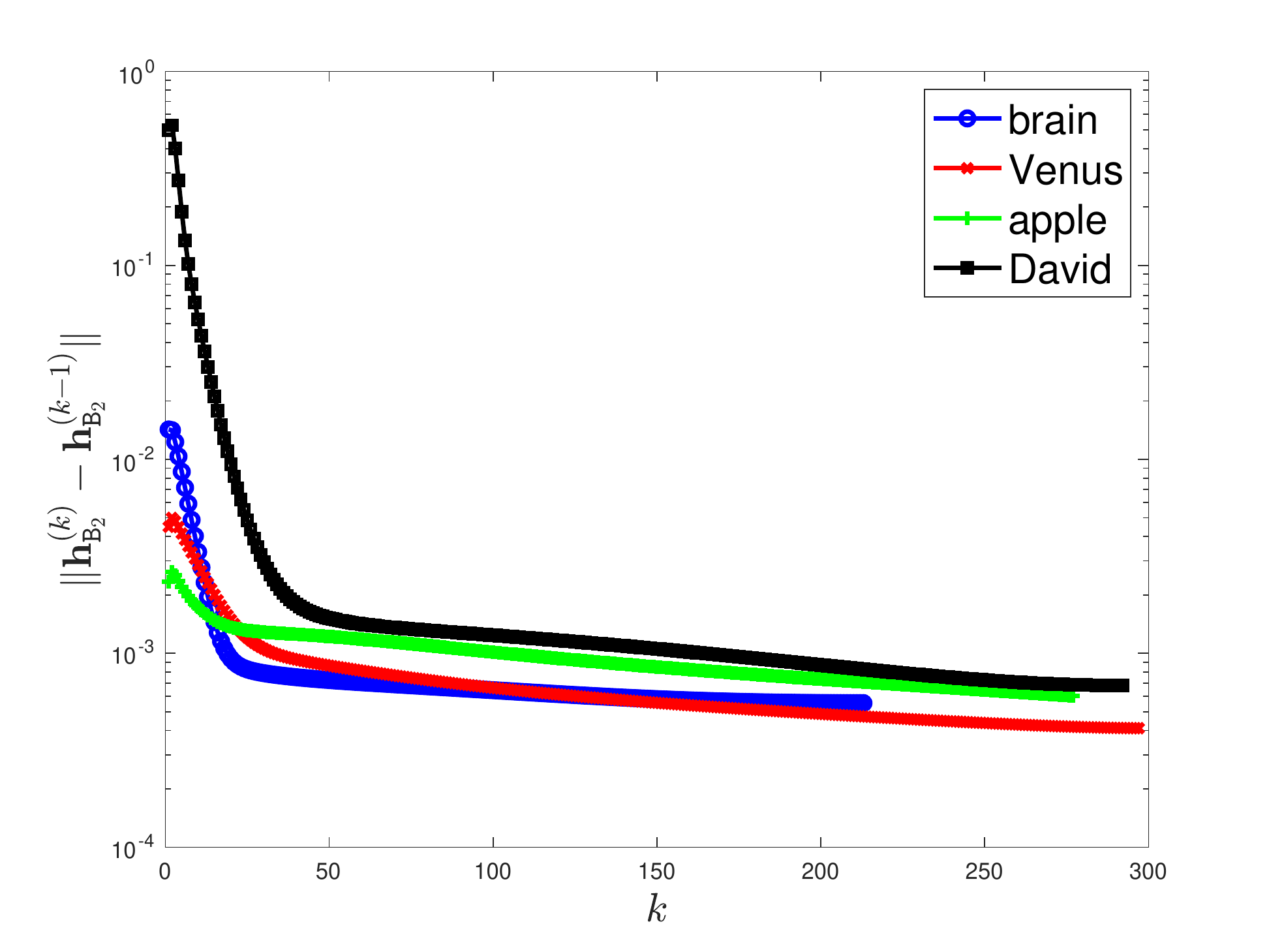}
\caption{$\| \h_{\B_2}^{(k)} - \h_{\B_2}^{(k-1)} \|_2$}
    \label{fig:h_B2_conv}
\end{subfigure}
\begin{subfigure}[b]{0.49\textwidth}
\center
    \includegraphics[width=\textwidth]{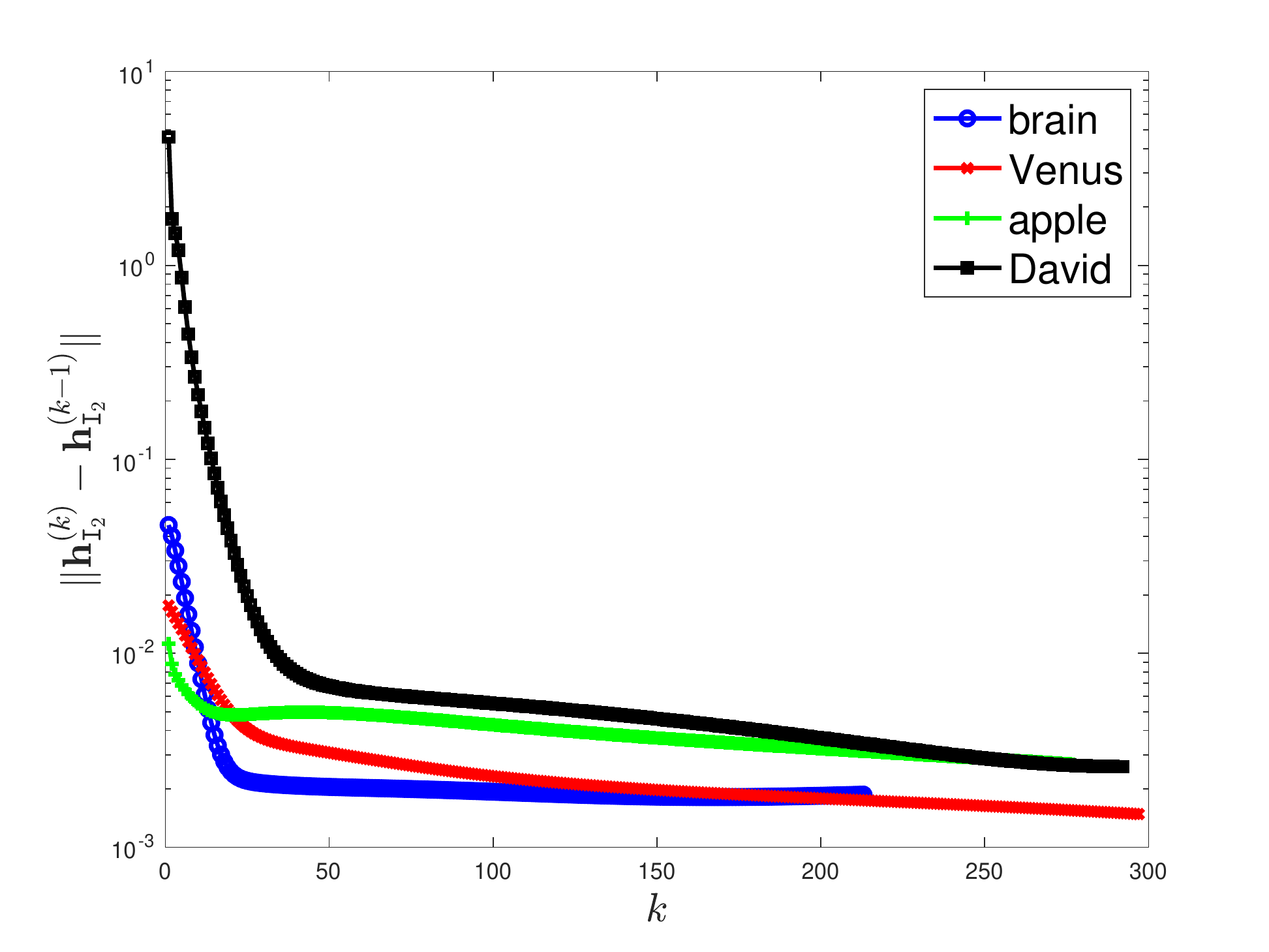}
\caption{$\| \h_{\I_2}^{(k)} - \h_{\I_2}^{(k-1)} \|_2$}
    \label{fig:h_I2_conv}
\end{subfigure}
\begin{subfigure}[b]{0.49\textwidth}
\center
    \includegraphics[width=\textwidth]{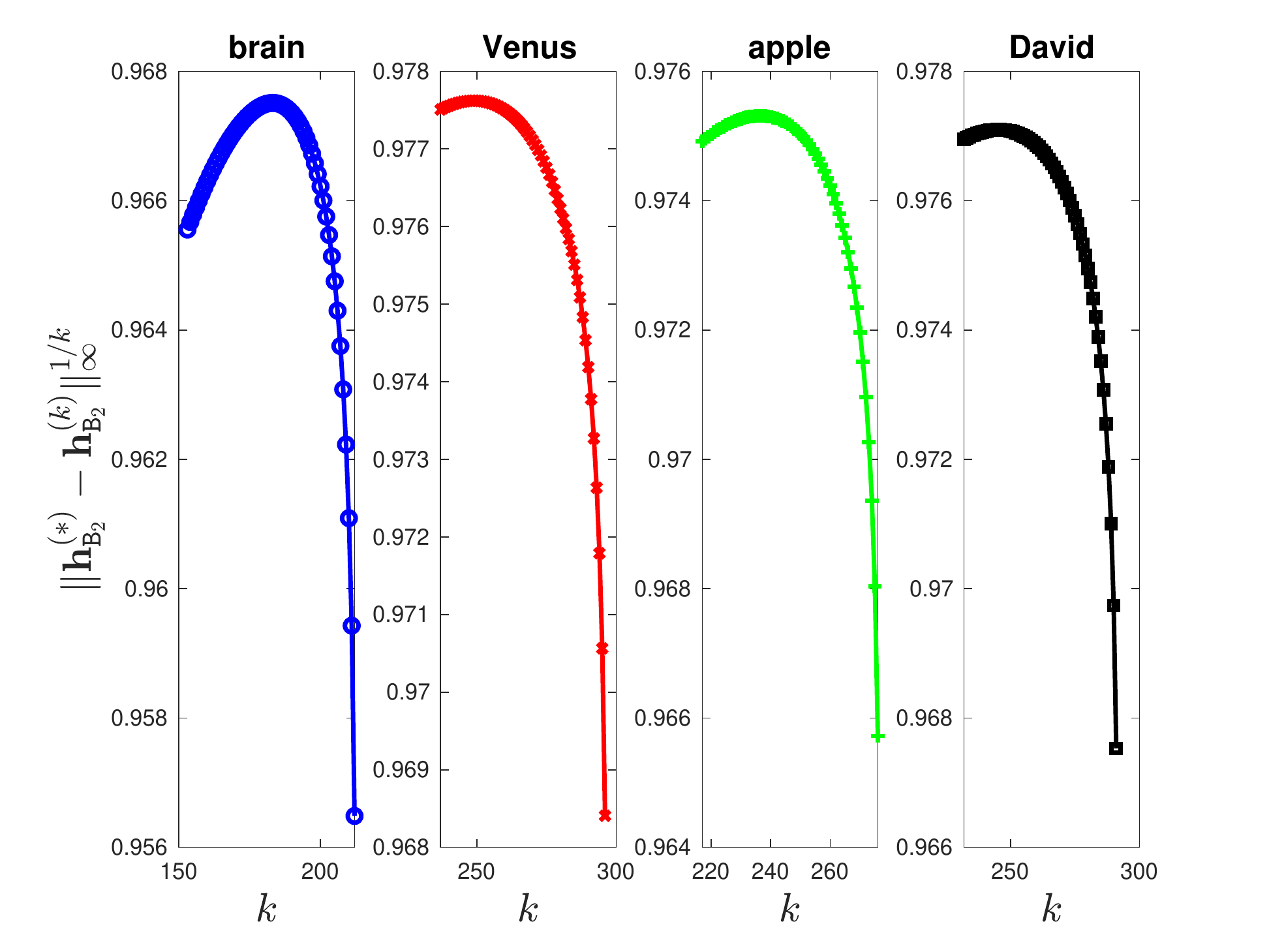}
\caption{$\| \h_{\B_2}^{(\ast)} - \h_{\B_2}^{(k)} \|_{\infty}^{1/k}$}
    \label{fig:h_B2_R_linear_conv}
\end{subfigure}
\begin{subfigure}[b]{0.49\textwidth}
\center
    \includegraphics[width=\textwidth]{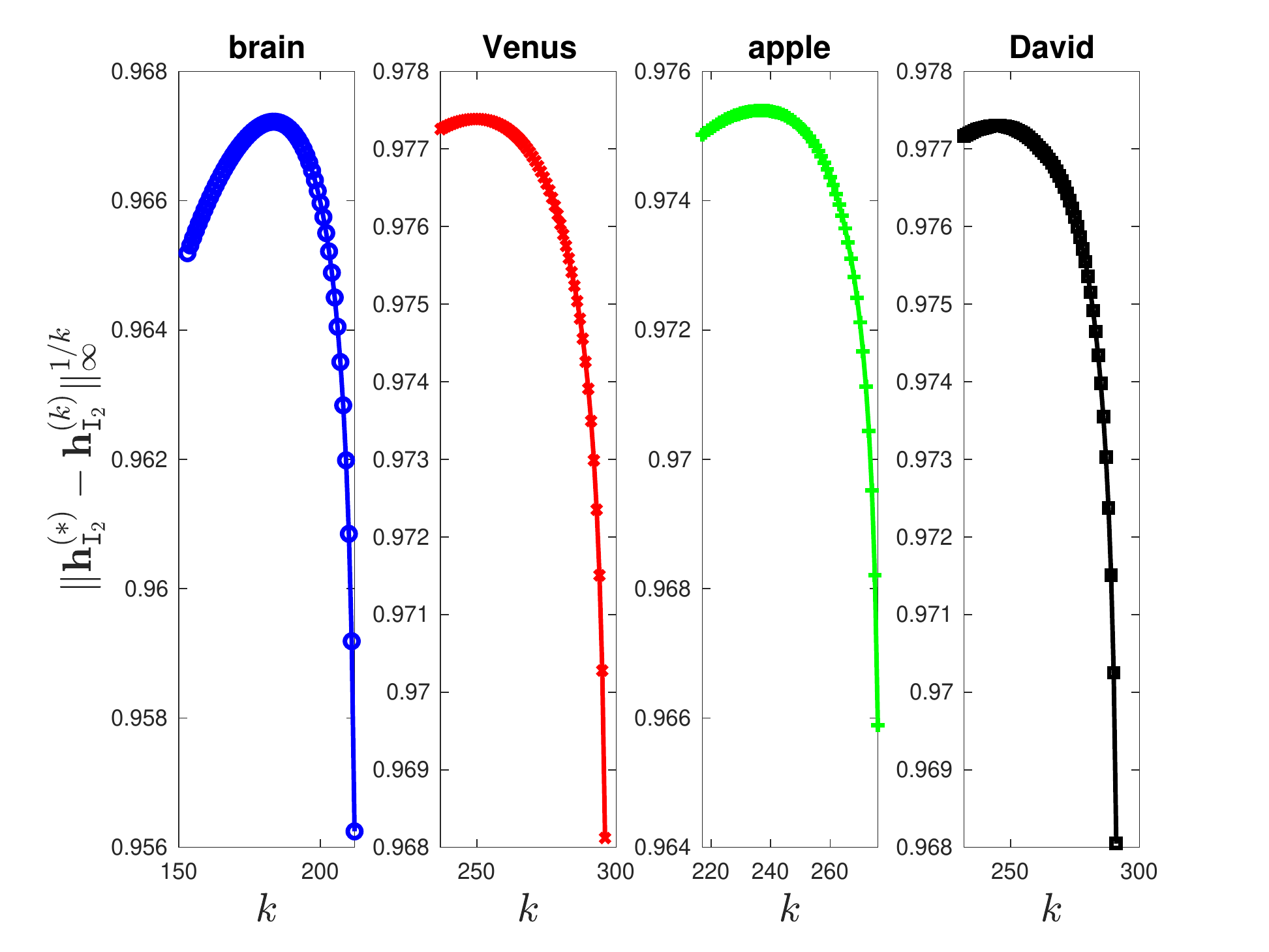}
\caption{$\| \h_{\I_2}^{(\ast)} - \h_{\I_2}^{(k)} \|_{\infty}^{1/k}$}
    \label{fig:h_I2_R_linear_conv}
\end{subfigure}
\caption{Convergence behavior of $\h_{\B_2}^{(k)}$ and $\h_{\I_2}^{(k)}$ in (A) and (B), and enlarged results of the R-linear convergences of $\h_{\B_2}^{(k)}$ and $\h_{\I_2}^{(k)}$ in (C) and (D).}
\label{fig:h_I2B2_R_linear_conv}
\end{figure}

In Figures \ref{fig:h_B2_conv} and \ref{fig:h_I2_conv}, we plot the results of $\| \h_{\B_2}^{(k)} - \h_{\B_2}^{(k-1)} \|_2$ and $\| \h_{\I_2}^{(k)} - \h_{\I_2}^{(k-1)} \|_2$, respectively, for each iteration $k$ to illustrate the convergence behavior of four benchmark models with ``mesh3''. The algorithm is stopped at $m$-th iteration for brain ($m = 213$), Venus ($m = 297$), apple ($m = 277$) and David ($m = 292$) models. Using $\h_{\B_2}^{(m)}$ and $\h_{\I_2}^{(m)}$ as the convergent vectors $\h_{\B_2}^{(\ast)}$ and $\h_{\I_2}^{(\ast)}$, respectively, we plot the results $\| \h_{\B_2}^{(\ast)} - \h_{\B_2}^{(k)} \|_{\infty}^{1/k}$ and $\| \h_{\I_2}^{(\ast)} - \h_{\I_2}^{(k)} \|_{\infty}^{1/k}$ with $k = m-60, \ldots, m-1$ in Figures~\ref{fig:h_B2_R_linear_conv} and \ref{fig:h_I2_R_linear_conv}. The results show that $\| \h_{\B_2}^{(\ast)} - \h_{\B_2}^{(k)} \|_{\infty}^{1/k}$ and $\| \h_{\I_2}^{(\ast)} - \h_{\I_2}^{(k)} \|_{\infty}^{1/k}$ increase to
nearly $1$ and then decreases linearly as $k$ increases, which indicates that the convergence of Algorithm~\ref{alg:SEM_fix_bnd_k} is asymptotically R-linear as shown the result (i) in Theorem~\ref{thm3.6}.

\section{Conclusions} \label{sec:5}
In this paper, we revisit the stretch energy minimization Algorithm~\ref{alg:SEM_fix_bnd_k} for the spherical equiareal parameterization between a simply connected closed surface $\mathcal{M}$ and $\mathbb{S}^2$. In order to compensate for deficiencies of convergence theory for Algorithm~\ref{alg:SEM_fix_bnd_k}, we intend to first provide a mathematical foundation in Theorem~\ref{invconformal} for the representation of the total area distortion (stretch energy) between $\mathcal{M}$ and $\mathbb{S}^2$ is equal to the sum of the Dirichlet energies on the unit discs corresponding to the southern  semi-sphere and the inverse of the northern hemi-sphere, respectively. In the light of this representation we are motivated to find the optimal solution for the stretch energy by minimizing the Dirichlet energies on the associated southern and northern discs alternatively, until convergence. Furthermore, we show that the determinant of Jacobian of the optimal solution is equal to one which is area-preserving that we expected. For Algorithm~\ref{alg:SEM_fix_bnd_k}, we first derive the complicated relation between errors at the $(k+1)$th and $k$th steps of Algorithm~\ref{alg:SEM_fix_bnd_k} which is a crucial expression for convergence. Then, under some mild and reasonable assumption, we prove Algorithm~\ref{alg:SEM_fix_bnd_k} either converges asymptotically and R-linearly, or forms a quasi-periodic solution. The mild assumption for convergence is guaranteed by numerical validations. Numerical experiments on various models from benchmarks demonstrate the asymptotically R-linear convergence and show the effectiveness and reliability of Algorithm~\ref{alg:SEM_fix_bnd_k}.

\section*{Acknowledgments}
This work was partially supported by the Ministry of Science and Technology (MoST), the National Center for Theoretical Sciences, and the ST Yau Center in Taiwan. W.-W. Lin and T.-M. Huang were partially supported by
MoST 110-2115-M-A49-004- and 110-2115-M-003-012-MY3, respectively.

\bibliographystyle{amsplain}
\bibliography{reference}
\end{document}